\documentclass[11pt]{article}
\usepackage[utf8]{inputenc}
\usepackage[T1]{fontenc}
\usepackage{lmodern}
\usepackage[a4paper,margin=25mm]{geometry}
\usepackage{amsmath,amssymb,amsfonts,amsthm,mathtools}
\usepackage{booktabs,graphicx,enumitem,adjustbox}
\usepackage[hidelinks]{hyperref}
\usepackage{microtype}
\newtheorem{theorem}{Theorem}[section]
\newtheorem{lemma}[theorem]{Lemma}
\newtheorem{proposition}[theorem]{Proposition}
\newtheorem{corollary}[theorem]{Corollary}
\theoremstyle{definition}
\newtheorem{definition}[theorem]{Definition}
\newtheorem{assumption}[theorem]{Assumption}
\theoremstyle{remark}
\newtheorem{remark}[theorem]{Remark}
\hypersetup{pdftitle={Second-order geometric consistency of nonlinear FE-HMM on surfaces via corrector non-degeneracy},pdfauthor={Marcial Nguemfouo},pdfkeywords={numerical homogenization, surface finite elements, geometric consistency}}

\title{Second-order geometric consistency of nonlinear FE--HMM on surfaces
via corrector non-degeneracy}
\author{Marcial Nguemfouo\thanks{Department of Mathematics and Applied
Mathematics, North-West University, Potchefstroom, South Africa, and
Department of Mathematics, University of Yaoundé I, Cameroon
(\href{mailto:marcial.nguemfouo@facsciences-uy1.cm}{\texttt{marcial.nguemfouo@facsciences-uy1.cm}}).}}
\date{November 5, 2026}

\begin{document}
\maketitle

\begin{abstract}
We analyze the geometric variational crime of nonlinear finite element
heterogeneous multiscale methods (FE--HMM) for $p$-growth monotone laws on
curved surfaces, $p\ge2$. The practical method projects discrete gradients
onto the exact tangent planes; we compare it with the exactly lifted method
on the same discrete space. Without further hypotheses, the difference is
$O(h_g^{p/(p-1)})$ in the metric induced by the cell correctors. The
optimal order $O(h_g^2)$ holds with a constant given by a weighted negative
moment of the corrected gradients, evaluated along the practical solution
only; no bounds on discrete gradients in $L^\infty$ are used. The negative
moment is a sufficient condition for the lower shifted $p$-structure of the
discrete homogenized flux, which for $p>2$ is not available in general. For
the model law and a fixed cell mesh, its uniform boundedness is equivalent
to a finite-dimensional directional criterion on degenerate correctors, and
for laminates it holds with an explicit bound, uniform in the cell mesh.
With macro--micro quasi-optimality and an estimate for coefficient
sampling, the reconstructed microscopic field converges as
$O(h+H^{\min\{s,1\}}+h_g^2)$, where $H^s$ is the cell approximation rate. Experiments on a torus with an exact closest-point lift
and a transported anisotropic $p=4$ material show second-order geometric
convergence up to $98\,304$ macroscopic quadrature points, each with a fully
solved nonlinear cell problem, with bounded a posteriori constants for both
load treatments.
\end{abstract}


\par\medskip

\section{Introduction}\label{sec:intro}

Heterogeneous multiscale finite element methods (FE--HMM) approximate a
homogenized problem without ever forming the effective law in closed form:
at every macroscopic quadrature point one solves a microscopic cell problem
and feeds the resulting effective flux into the macroscopic residual
\cite{EEngquist03,AbdulleEEngquistVandenEijnden12}. For nonlinear monotone
laws this approach is well understood in Euclidean domains
\cite{Gloria06,AbdulleVilmart14,AbdulleHuber16,HenningOhlberger15,Verfurth22}.
On curved surfaces a third discretization enters: the surface itself is
replaced by a polyhedral approximation $\Gamma_h$, which produces a
geometric variational crime \cite{Dziuk88,Demlow09,DziukElliott13}. For
linear diffusion this crime is of second order in the geometric mesh size,
$O(h^2)$, and this order is known to be sharp.

This paper answers the corresponding question for nonlinear FE--HMM with
$p$-growth on surfaces: \emph{when is the geometric error of second order,
and what structural property of the discrete homogenized law decides it?}
The question is not a formality. The natural error measure for $p$-growth
problems is quadratic in the nonlinear variable $F_\delta(\xi)=(\delta+|\xi|)^{(p-2)/2}\xi$
\cite{BarrettLiu94,DieningEttwein08,DieningKreuzer08}. An energy comparison
of the exact and approximate surface functionals transfers the $O(h^2)$
metric perturbation into the squared natural distance and therefore yields
only $O(h)$. Recovering $O(h^2)$ requires an operator-level argument, and
such an argument needs two-sided bounds on the effective flux in the shifted
form
\begin{equation}\label{eq:shifted-intro}
\begin{aligned}
\big(q_H(\xi)-q_H(\eta)\big)\cdot(\xi-\eta)&\gtrsim(\delta+|\xi|+|\eta|)^{p-2}|\xi-\eta|^2,
\\
|q_H(\xi)-q_H(\eta)|&\lesssim(\delta+|\xi|+|\eta|)^{p-2}|\xi-\eta|,
\end{aligned}
\end{equation}
uniformly in the cell discretization. The microscopic law satisfies
\eqref{eq:shifted-intro} by assumption; whether the \emph{discrete
homogenized} law $q_H$ does is a genuine question about nonlinear
homogenization. The transfer of growth and monotonicity conditions to the
homogenized operator goes back to \cite{FuscoMoscariello87,ChiadoPiatDefranceschi90};
what the surface analysis needs is the shifted form, with constants
independent of $\delta$ and of the cell mesh.

\subsection*{Contributions}
\begin{enumerate}[label=(\roman*),leftmargin=2em]
\item \textbf{Second-order geometric consistency with an a posteriori
constant} (Theorem~\ref{thm:geom}). This is the main result. Comparing the
practical surface method, which projects discrete gradients, with the
exactly lifted method on the same discrete space, we prove an
\emph{unconditional} bound of order $h_g^{p/(p-1)}$ in the corrector-induced
metric and the \emph{optimal} order $h_g^2$ with a constant given by a
weighted negative moment of the corrected gradients evaluated along the
practical solution only. Only $W^{1,p}$ stability of the discrete solutions
is used: no $L^\infty$ bounds on discrete gradients and no regularity of
the corrector map.

\item \textbf{Lower shifted structure via a negative moment of the
correctors} (Theorem~\ref{thm:lower}, Proposition~\ref{prop:reduction}). For
$p>2$ the lower bound in \eqref{eq:shifted-intro} is not known to hold in
general and is assumed in the homogenization literature
\cite[Remark~2.5]{WangXuZhao19}. We show that it holds with a constant
controlled by a single-state negative moment $\Lambda_H$ of the corrected
gradients; this is a sufficient condition, weaker than the pointwise
non-degeneracy assumed in \cite{Wolf22}. Under it the corrector-induced
metric, in which macro--micro methods are quasi-optimal, is equivalent to
the classical natural distance (Corollary~\ref{cor:metric}). For the model
law and a fixed cell mesh, boundedness of $\Lambda_H$ uniformly in $\delta$
and $\xi$ is equivalent to a finite-dimensional directional criterion on the
degenerate correctors (Theorem~\ref{thm:directions}), which can be
certified with finitely many cell solves (Lemma~\ref{lem:certify}). For
laminates we give an explicit bound, uniform in $\delta$ and in the cell
mesh (Proposition~\ref{prop:laminate}). The upper half of
\eqref{eq:shifted-intro} is known for the continuous operator
\cite{WangXuZhao19,ClozeauGloria23}; we record its discrete, $\delta$-uniform
version (Theorem~\ref{thm:upper}).

\item \textbf{Complete estimate} (Corollary~\ref{cor:total}). With the
macro--micro quasi-optimality, proved in detail on the surface, the
reconstructed microscopic field converges as $h+H^{\min\{s,1\}}+h_g^2$, where
$H^s$ is the cell approximation rate, including the
error from sampling the coefficient on the cell mesh
(Lemma~\ref{lem:sampling}); under a continuous index bound, so does the
classical natural distance of the macroscopic gradients.

\item \textbf{Numerics} (Section~\ref{sec:num}). On a torus with the
closed-form closest-point lift and a transported anisotropic $p=4$
material, the geometric error converges with experimental orders between
$1.97$ and $2.00$ up to $98\,304$ macroscopic quadrature points, each with a
fully solved nonlinear cell problem, with a bounded a posteriori constant,
for both load treatments and $\delta$ down to $0.02$.
\end{enumerate}

The structural results of Section~\ref{sec:structure} are not specific to
surfaces. Any a priori analysis of a nonlinear multiscale method that uses
the discrete effective flux $q_H$ in the natural distance, for instance
under quadrature, coefficient perturbation or model reduction of the cell
problems, needs the lower shifted bound for $q_H$; the negative moment
$\Lambda_H$ and the directional criterion provide a checkable route to it.

\subsection*{Relation to previous work}
The FE--HMM error analysis for nonlinear monotone problems
\cite{AbdulleVilmart14,AbdulleHuber16}, the variational framework of
\cite{Gloria06} and a posteriori control \cite{HenningOhlberger15} are
Euclidean. Structure conditions for homogenized monotone operators go back
to \cite{FuscoMoscariello87,ChiadoPiatDefranceschi90}; the shifted upper
bound for the continuous operator is stated in
\cite[Lemma~2.4]{WangXuZhao19} and \cite{ClozeauGloria23}, while the shifted
lower bound for $p>2$ is imposed as an assumption in
\cite{WangXuZhao19}. For the $p$-Laplacian with a local defect, Wolf
\cite{Wolf22} assumes the pointwise non-degeneracy
$|\xi+\nabla w_\xi^{\rm per}|\ge c|\xi|$ of the continuous periodic correctors,
which is our condition \eqref{eq:kappa} for $W_0$, proves it for laminates,
and discusses a weaker weighted Poincaré condition implied by
$|\xi+\nabla w^{\rm per}_\xi|^{2-p}\in L^{d/2}$. Our averaged index
$\Lambda_0$ is the $L^1$ version of this negative moment, the two coincide
for $d=2$, and we use it for a different purpose: as a constant in a priori
estimates for the discrete method. HMM on rough surfaces \cite{AbdulleSchwab05} treats linear
diffusion after transformation to a parameter domain. Geometric
homogenization theory on manifolds
\cite{Pak05,DobberschutzBohm12,Dobberschutz18} provides the continuous
limit but no discretization. The macro--micro energy arguments we use in
Section~\ref{sec:setting} follow from cell stationarity and conformity and
are not claimed as new. The new elements are the results of Section~\ref{sec:structure} on the
lower shifted structure, which apply to any nonlinear FE--HMM, and their
use to obtain the optimal geometric order with stability-only assumptions.

\subsection*{Outline}
Section~\ref{sec:setting} fixes the microscopic structure and recalls the
discrete homogenized law. Section~\ref{sec:structure} proves the inheritance
results. Section~\ref{sec:geometry} contains the surface analysis.
Section~\ref{sec:num} reports the numerical experiments, and
Section~\ref{sec:conclusion} concludes.

\section{Setting and the discrete homogenized law}\label{sec:setting}

\subsection{Shifted \texorpdfstring{$p$}{p}-structure}
Let $p\ge2$, $\delta\ge0$, $m\ge1$ and $Y=(0,1)^m$ with periodic
identification. For $t\ge0$ and $\xi,\eta\in\mathbb R^m$ write
\[
W(\xi,\eta):=\delta+|\xi|+|\eta|,\qquad
F_\delta(\xi):=(\delta+|\xi|)^{\frac{p-2}{2}}\xi .
\]

\begin{assumption}[Microscopic law]\label{ass:micro}
The map $a:\Omega\times Y\times\mathbb R^m\to\mathbb R^m$ is Carathéodory,
$Y$-periodic in $y$, $a(x,y,0)=0$, and there are constants
$0<c_1\le C_1$ such that for a.e.\ $(x,y)$ and all $z,w\in\mathbb R^m$
\begin{align}
\big(a(x,y,z)-a(x,y,w)\big)\cdot(z-w)&\ge c_1\,W(z,w)^{p-2}|z-w|^2,\label{eq:mono}\\
|a(x,y,z)-a(x,y,w)|&\le C_1\,W(z,w)^{p-2}|z-w|.\label{eq:lip}
\end{align}
\end{assumption}

These are the standard shifted conditions for $p$-Laplace-type operators;
they imply $|F_\delta(z)-F_\delta(w)|^2\simeq W(z,w)^{p-2}|z-w|^2$ with
constants depending only on $p$ \cite{DieningEttwein08,DieningKreuzer08}.
The model of our experiments,
\begin{equation}\label{eq:model}
a(x,y,z)=\gamma(x)\big(\delta^2+z^\top A(y)z\big)^{\frac{p-2}{2}}A(y)z,\qquad
\lambda_-I\preceq A(y)\preceq\lambda_+I,\quad \gamma_-\le\gamma\le\gamma_+,
\end{equation}
satisfies Assumption~\ref{ass:micro}; it is the gradient of
$\phi(x,y,z)=\frac{\gamma(x)}{p}(\delta^2+z^\top A(y)z)^{p/2}$. In this case the
microscopic natural variable is
$F_A(x,y,z)=\gamma(x)^{1/2}(\delta^2+z^\top A(y)z)^{\frac{p-2}{4}}A(y)^{1/2}z$, and
$|F_A(z)-F_A(w)|^2\simeq W(z,w)^{p-2}|z-w|^2$.

\subsection{Discrete cell problems}
Let $W_H\subset W^{1,p}_{\#}(Y)$ be any conforming space of continuous
periodic finite element functions with zero mean, and set
$W_0=W^{1,p}_\#(Y)$ for the continuous problem. For $\xi\in\mathbb R^m$,
$\pi_{\xi,H}\in W_H$ solves
\begin{equation}\label{eq:cell}
\int_Y a\big(x,y,\xi+\nabla_y\pi_{\xi,H}\big)\cdot\nabla_y\vartheta\,dy=0
\qquad\forall\vartheta\in W_H,
\end{equation}
and we write $Z_H(x,y,\xi)=\xi+\nabla_y\pi_{\xi,H}$ and
$q_H(x,\xi)=\int_Y a(x,y,Z_H(x,y,\xi))\,dy$. Existence and uniqueness follow
from the Browder--Minty theorem. The continuous corrector $Z$ and flux $q$
correspond to $W_0$. Testing the two cell equations with
$\pi_{\xi,H}-\pi_{\eta,H}\in W_H$ gives the basic identity
\begin{equation}\label{eq:basic}
\big(q_H(\xi)-q_H(\eta)\big)\cdot(\xi-\eta)
=\int_Y\big(a(Z_H(\xi))-a(Z_H(\eta))\big)\cdot\big(Z_H(\xi)-Z_H(\eta)\big)\,dy ,
\end{equation}
where here and below the arguments $x,y$ are omitted. Since
$\int_Y(Z_H(\xi)-Z_H(\eta))\,dy=\xi-\eta$, Jensen's inequality and
\eqref{eq:mono} give $(q_H(\xi)-q_H(\eta))\cdot(\xi-\eta)\ge c_1|\xi-\eta|^p$.

The \emph{corrector-induced metric} is
\begin{equation}\label{eq:dxH}
d_{x,H}(\xi,\eta)^2:=\int_Y\big|F_\delta(Z_H(\xi))-F_\delta(Z_H(\eta))\big|^2dy ,
\end{equation}
with $F_\delta$ replaced by $F_A$ for \eqref{eq:model} (the constants in
the equivalences below then also depend on $\lambda_\pm,\gamma_\pm$). By \eqref{eq:basic},
\eqref{eq:mono}, \eqref{eq:lip} and the equivalence above,
\begin{equation}\label{eq:dxH-equiv}
d_{x,H}(\xi,\eta)^2\simeq\int_Y W(Z_H(\xi),Z_H(\eta))^{p-2}|Z_H(\xi)-Z_H(\eta)|^2dy
\simeq\big(q_H(\xi)-q_H(\eta)\big)\cdot(\xi-\eta).
\end{equation}
In the variational case the right-hand side is the symmetrized Bregman
distance of the discrete effective energy, and $d_{x,H}$ is the metric in
which the macro--micro method is quasi-optimal; see
\cite{Gloria06,AbdulleHuber16} and Section~\ref{sec:geometry}.

\subsection{Shifted stability of discrete correctors}

\begin{lemma}\label{lem:stab}
There is $C_s=C_s(p,c_1,C_1)$, independent of $\delta$, $H$, $x$ and $\xi$,
such that
\[
\int_Y\big(\delta+|Z_H(\xi)|\big)^p\,dy\le C_s\,(\delta+|\xi|)^p .
\]
\end{lemma}

\begin{proof}
Write $Z=Z_H(\xi)$ and $\psi(t)=(\delta+t)^{p-2}t^2$. Testing \eqref{eq:cell}
with $\pi_{\xi,H}$ gives $\int_Y a(Z)\cdot Z=\int_Y a(Z)\cdot\xi$. With $w=0$,
\eqref{eq:mono} and \eqref{eq:lip} give $a(Z)\cdot Z\ge c_1\psi(|Z|)$ and
$|a(Z)|\le C_1(\delta+|Z|)^{p-2}|Z|$. Hence
\[
c_1\int_Y\psi(|Z|)\le C_1|\xi|\int_Y(\delta+|Z|)^{p-2}|Z|
\le \tfrac{c_1}{2}\int_Y\psi(|Z|)+C\int_Y(\delta+|Z|)^{p-2}|\xi|^2 .
\]
Since $(\delta+|Z|)^{p-2}\le 2^{p-3}\big((\delta+|\xi|)^{p-2}+|Z|^{p-2}\big)$ for $p\ge3$
(and trivially for $2\le p<3$ with a modified constant), and
$|Z|^{p-2}|\xi|^2\le\varepsilon|Z|^p+C_\varepsilon|\xi|^p\le\varepsilon\psi(|Z|)+C_\varepsilon|\xi|^p$
by Young's inequality with exponents $\frac{p}{p-2},\frac p2$, absorption
gives $\int_Y\psi(|Z|)\le C(\delta+|\xi|)^{p-2}|\xi|^2$. Finally
$(\delta+t)^p\le2^{p-1}(\delta^p+t^p)\le2^{p-1}(\delta^p+\psi(t))$.
\end{proof}

\section{Shifted structure of the discrete homogenized law}\label{sec:structure}

\begin{theorem}[Upper shifted structure is inherited]\label{thm:upper}
Under Assumption~\ref{ass:micro} there is $C_{\rm up}=C_{\rm up}(p,c_1,C_1)$,
independent of $\delta$, the cell space $W_H$ and $x$, such that for all
$\xi,\eta\in\mathbb R^m$
\[
|q_H(\xi)-q_H(\eta)|\le C_{\rm up}\,W(\xi,\eta)^{p-2}|\xi-\eta| .
\]
The same holds for the continuous homogenized flux $q$. In particular
$|q_H(\xi)|\le C_{\rm up}(\delta+|\xi|)^{p-2}|\xi|$.
\end{theorem}

For the continuous flux $q$ this bound is known
\cite[Lemma~2.4]{WangXuZhao19}, \cite{ClozeauGloria23}. We include the
short argument because the surface analysis needs it uniformly in the cell
space and in $\delta$.

\begin{proof}
Let $Z_\xi=Z_H(\xi)$, $Z_\eta=Z_H(\eta)$, $D=Z_\xi-Z_\eta$,
$w=W(Z_\xi,Z_\eta)$ and $X=|q_H(\xi)-q_H(\eta)|$. By \eqref{eq:lip} and the
Cauchy--Schwarz inequality,
\[
X\le C_1\int_Y w^{p-2}|D|\le C_1\Big(\int_Y w^{p-2}\Big)^{1/2}\Big(\int_Y w^{p-2}|D|^2\Big)^{1/2}.
\]
By Hölder's inequality and Lemma~\ref{lem:stab},
$\int_Y w^{p-2}\le(\int_Y w^p)^{(p-2)/p}\le C\,W(\xi,\eta)^{p-2}$. By
\eqref{eq:mono} and \eqref{eq:basic},
$c_1\int_Y w^{p-2}|D|^2\le(q_H(\xi)-q_H(\eta))\cdot(\xi-\eta)\le X|\xi-\eta|$.
Therefore $X\le C\,W(\xi,\eta)^{(p-2)/2}(X|\xi-\eta|)^{1/2}$, that is
$X\le C\,W(\xi,\eta)^{p-2}|\xi-\eta|$. The bound for $q_H(\xi)$ follows from
$q_H(0)=0$, which holds because $a(\cdot,\cdot,0)=0$ makes $\pi_{0,H}=0$
the cell solution.
\end{proof}

The argument uses the flux defect on both sides: the weighted energy of the
corrector difference is controlled by the flux defect through
\eqref{eq:basic}, which closes the estimate. For the lower bound no such
closure is available, because a weighted average of $|D|^2$ can be small
while $|\int_Y D|$ is not. For $p>2$ the lower shifted bound for the
homogenized operator is, to our knowledge, not available in general and is
imposed as a hypothesis in \cite[Remark~2.5]{WangXuZhao19}. This leads to
the following index.

\begin{definition}[Non-degeneracy index]\label{def:theta}
For $\xi\ne\eta$ set
\[
\Theta_H(x;\xi,\eta):=W(\xi,\eta)^{p-2}\int_Y W\big(Z_H(x,y,\xi),Z_H(x,y,\eta)\big)^{2-p}dy,
\]
\[
\Theta_H^*:=\operatorname*{ess\,sup}_{x}\ \sup_{\xi\ne\eta}\Theta_H(x;\xi,\eta).
\]
For $p=2$, $\Theta_H\equiv1$.
\end{definition}

\begin{theorem}[Lower shifted structure]\label{thm:lower}
Under Assumption~\ref{ass:micro}, for all $x$ and $\xi\ne\eta$,
\[
\big(q_H(\xi)-q_H(\eta)\big)\cdot(\xi-\eta)\ \ge\ \frac{c_1}{\Theta_H(x;\xi,\eta)}\,W(\xi,\eta)^{p-2}|\xi-\eta|^2 .
\]
If the corrected gradients satisfy the pointwise non-degeneracy condition
\begin{equation}\label{eq:kappa}
|Z_H(x,y,\xi)|\ge\kappa\,|\xi|\qquad\text{for a.e. }y\in Y,\ \text{all }\xi,
\end{equation}
for some $\kappa\in(0,1]$, then $\Theta_H^*\le\kappa^{2-p}$. The same
statements hold for the continuous correctors $Z$, with index $\Theta_0$.
\end{theorem}

\begin{proof}
With the notation of the previous proof, Cauchy--Schwarz gives
\[
|\xi-\eta|^2=\Big|\int_Y D\Big|^2\le\int_Y w^{p-2}|D|^2\int_Y w^{2-p}.
\]
The second
factor equals $\Theta_H W(\xi,\eta)^{2-p}$, and \eqref{eq:mono},
\eqref{eq:basic} bound the first one by $c_1^{-1}(q_H(\xi)-q_H(\eta))\cdot(\xi-\eta)$.
Under \eqref{eq:kappa},
$w\ge\delta+\kappa(|\xi|+|\eta|)\ge\kappa W(\xi,\eta)$ pointwise, so
$\int_Y w^{2-p}\le\kappa^{2-p}W(\xi,\eta)^{2-p}$ because $2-p\le0$.
\end{proof}

\begin{corollary}[Equivalence of natural metrics]\label{cor:metric}
Under Assumption~\ref{ass:micro},
\[
c\,\Theta_H(x;\xi,\eta)^{-1}|F_\delta(\xi)-F_\delta(\eta)|^2\ \le\ d_{x,H}(\xi,\eta)^2\ \le\ C\,|F_\delta(\xi)-F_\delta(\eta)|^2 ,
\]
with $c,C$ depending only on $p,c_1,C_1$. In particular, the
corrector-induced metric is always dominated by the classical natural
distance, and it is equivalent to it whenever $\Theta_H^*<\infty$.
\end{corollary}

\begin{proof}
Combine \eqref{eq:dxH-equiv} with Theorems~\ref{thm:upper},
\ref{thm:lower}, the bound $(q_H(\xi)-q_H(\eta))\cdot(\xi-\eta)\le
|q_H(\xi)-q_H(\eta)||\xi-\eta|$ and the equivalence
$|F_\delta(\xi)-F_\delta(\eta)|^2\simeq W(\xi,\eta)^{p-2}|\xi-\eta|^2$.
\end{proof}

\begin{remark}[Computability]\label{rem:computable}
In an FE--HMM computation the corrected gradients $Z_H(x_j,\cdot,\xi)$ are
available at every macroscopic quadrature point $x_j$ and every macroscopic
state $\xi$ visited by the solver, so every index introduced below can be
evaluated along the computation without additional cell solves.\end{remark}

\subsection{A single-state index and an exact finite-dimensional test}\label{sec:single}
The pair index $\Theta_H$ can be replaced by a quantity that involves one
macroscopic state and only an integral of the corrected gradient.

\begin{definition}[Single-state index]\label{def:lambda}
For $\xi\in\mathbb R^m$ set
\[
\Lambda_H(x;\xi):=(\delta+|\xi|)^{p-2}\int_Y\big(\delta+|Z_H(x,y,\xi)|\big)^{2-p}dy\ \in(0,\infty],
\]
\[
\Lambda_H^*:=\operatorname*{ess\,sup}_x\sup_{\xi}\Lambda_H(x;\xi),
\]
for all states with $\delta+|\xi|>0$, and $\Lambda_0$, $\Lambda_0^*$ for the
continuous correctors. For $\delta=0$ the state $\xi=0$ is excluded from the
supremum; wherever the weighted quantity $\Lambda_H(x;\xi)(\delta+|\xi|)^p$
appears, it is set to $0$ at $\delta+|\xi|=0$, which is its limit when
$\Lambda_H$ is bounded. For $p=2$, $\Lambda_H\equiv1$.
\end{definition}

\begin{proposition}[Reduction]\label{prop:reduction}
Under Assumption~\ref{ass:micro}, for all $x$ and $\xi\ne\eta$ with
$|\xi|\ge|\eta|$,
\[
\Theta_H(x;\xi,\eta)\le2^{p-2}\Lambda_H(x;\xi).
\]
Moreover $\Lambda_H(x;\xi)\le\kappa^{2-p}$ under \eqref{eq:kappa}, and
$\Lambda_H(x;\xi)\le(1+|\xi|/\delta)^{p-2}$ for every $\delta>0$.
\end{proposition}

\begin{proof}
Since $p\ge2$ and $W(Z_H(\xi),Z_H(\eta))\ge\delta+|Z_H(\xi)|$,
\[\int_YW(Z_H(\xi),Z_H(\eta))^{2-p}\le\int_Y(\delta+|Z_H(\xi)|)^{2-p},\] and
$W(\xi,\eta)\le2(\delta+|\xi|)$. The second claim follows from
$\delta+|Z_H|\ge\delta+\kappa|\xi|\ge\kappa(\delta+|\xi|)$, the third from
$\delta+|Z_H|\ge\delta$.
\end{proof}

The last bound shows that non-degeneracy is only an issue in the regime
$|\xi|\gg\delta$, where the law is close to its degenerate limit. For the
model law this regime reduces to a compact set of directions.

\begin{theorem}[Uniform index by a directional criterion]\label{thm:directions}
Let $a$ be the model law \eqref{eq:model}, with $A=A(y)$, and let $W_H$ be
the $P_1$ space on a periodic simplicial mesh $\mathcal T_H$ of $Y$. For
$s\ge0$ and $\omega\in S^{m-1}$ let $Z_H^{s}(\cdot,\omega)$ denote the
discrete corrected gradient of the law with regularization $s$ at the
macroscopic state $\omega$; for $s=0$ this is the degenerate law
$(z^\top Az)^{\frac{p-2}{2}}Az$. Then:
\begin{enumerate}[label=(\roman*),leftmargin=2em]
\item (Homogeneity.) $Z_H(y,\xi)=|\xi|\,Z_H^{\delta/|\xi|}(y,\xi/|\xi|)$ and
$\Lambda_H(\xi)=\lambda_H(\delta/|\xi|,\xi/|\xi|)$ with
$\lambda_H(s,\omega):=(1+s)^{p-2}\int_Y(s+|Z_H^s(y,\omega)|)^{2-p}dy$.
\item (Characterization.) For $p>2$, $\sup\{\Lambda_H(\xi):\delta\ge0,\ \xi\ne0\}<\infty$
if and only if the finite-dimensional directional criterion
\begin{equation}\label{eq:test}
\kappa_H^0:=\min_{\omega\in S^{m-1}}\ \min_{T\in\mathcal T_H}|Z_H^0|_T(\omega)|>0 .
\end{equation}
holds. In that case there exists $s_1>0$ such that
$g_1:=\min\{|Z_H^s|_T(\omega)|:\,s\le s_1,\ \omega\in S^{m-1},\ T\in\mathcal T_H\}>0$,
and for every such $s_1$,
$\Lambda_H^*\le\max\big\{(1+s_1)^{p-2}g_1^{2-p},(1+s_1^{-1})^{p-2}\big\}$.
\end{enumerate}
Since $Z_H^0(\cdot,-\omega)=-Z_H^0(\cdot,\omega)$, for $m=2$ the criterion
\eqref{eq:test} involves a half-circle of directions. It does not involve
$\delta$, the macroscopic problem or the geometry. The minimum over the
continuum of directions is not a finite computation; Lemma~\ref{lem:certify}
turns finitely many cell solves into a certificate.
\end{theorem}

\begin{proof}
(i) The discrete cell solution minimizes
$w\mapsto\int_Y\frac1p(\delta^2+(\xi+\nabla w)^\top A(\xi+\nabla w))^{p/2}$ over
$W_H$. For $t>0$, replacing $(\delta,\xi,w)$ by $(t\delta,t\xi,tw)$
multiplies the functional by $t^p$, so minimizers correspond and
$Z_H$ is positively one-homogeneous in $(\delta,\xi)$. Take $t=|\xi|^{-1}$.

(ii) On $[0,\infty)\times S^{m-1}\times W_H$ the discrete energy is jointly
continuous, strictly convex and uniformly coercive in $w$ for $(s,\omega)$ in
compact sets; this holds also for $s=0$, because $z\mapsto|A^{1/2}z|^p$ is
strictly convex for $p>1$. Hence the minimizer, and with it each element
value $Z_H^s|_T(\omega)$, depends continuously on $(s,\omega)$. If
\eqref{eq:test} holds, the continuous function
$g(s,\omega)=\min_T|Z_H^s|_T(\omega)|$ is positive on $\{0\}\times S^{m-1}$,
hence on $[0,s_1]\times S^{m-1}$ for some $s_1>0$ by compactness, with
minimum $g_1>0$. For $s\le s_1$, $\lambda_H\le(1+s_1)^{p-2}g_1^{2-p}$; for
$s\ge s_1$, $\lambda_H\le((1+s)/s)^{p-2}\le(1+s_1^{-1})^{p-2}$. Conversely, if
$Z_H^0|_T(\omega_0)=0$ for some $T$, then $|Z_H^s|_T(\omega_0)|\to0$ as
$s\to0$ and $\lambda_H(s,\omega_0)\ge|T|(s+|Z_H^s|_T(\omega_0)|)^{2-p}\to\infty$.
Oddness follows from uniqueness, because the law is odd in $z$.
\end{proof}

\begin{lemma}[Certification between directions]\label{lem:certify}
Under the assumptions of Theorem~\ref{thm:directions}, let $c_1,C_{\rm up}$
be the constants of \eqref{eq:mono} and Theorem~\ref{thm:upper} for the
degenerate law, and let $\omega_k\in S^{m-1}$ with
$g_k:=\min_T|Z_H^0|_T(\omega_k)|>0$. Then for every $\omega\in S^{m-1}$ and
$T\in\mathcal T_H$,
\[
|Z_H^0|_T(\omega)|\ \ge\ g_k-\Big(\frac{2^{p-2}C_{\rm up}}{c_1\,|T|}\Big)^{1/2}g_k^{-\frac{p-2}{2}}\,|\omega-\omega_k| .
\]
Consequently, if the balls of radius $\rho_k$ around directions
$\omega_1,\dots,\omega_n$ cover the half-sphere and
\[
\min_k\Big(g_k-\Big(\frac{2^{p-2}C_{\rm up}}{c_1\min_T|T|}\Big)^{1/2}g_k^{-\frac{p-2}{2}}\rho_k\Big)>0,
\]
then \eqref{eq:test} holds. The lemma concerns the exact solutions of the
discrete cell problems; a rigorous numerical certificate must also account
for the error of the nonlinear solves.
\end{lemma}

\begin{proof}
Let $D=Z_H^0(\omega)-Z_H^0(\omega_k)$ and $w=|Z_H^0(\omega)|+|Z_H^0(\omega_k)|\ge g_k$
elementwise. By \eqref{eq:mono}, \eqref{eq:basic} and
Theorem~\ref{thm:upper} with $\delta=0$ and $|\omega|=|\omega_k|=1$,
$c_1|T|\,w_T^{p-2}|D_T|^2\le c_1\int_Yw^{p-2}|D|^2\le(q_H(\omega)-q_H(\omega_k))\cdot(\omega-\omega_k)\le2^{p-2}C_{\rm up}|\omega-\omega_k|^2$.
Hence $|D_T|\le(2^{p-2}C_{\rm up}/(c_1|T|))^{1/2}g_k^{-(p-2)/2}|\omega-\omega_k|$ and the
claim follows from the triangle inequality.
\end{proof}

A computer-assisted certificate would also require explicit values of
$c_1$ and $C_{\rm up}$. The certificate is conservative because of the
factor $|T|^{-1/2}$. We do not claim a numerical certification: in
Section~\ref{sec:num-directions} we report sampled minima, which detect
degeneracy but do not certify its absence.

\begin{remark}[What remains open]\label{rem:open}
Theorem~\ref{thm:directions} gives, for each cell mesh, an exact
criterion for an index bound that is uniform in $\delta$ and $\xi$, but the
bound may depend on $H$. A bound uniform in $H$ would follow from a
quantitative lower bound on the continuous degenerate correctors together
with $W^{1,\infty}$-convergence of the discrete ones. For linear planar
conductivity the corrector map is univalent with almost everywhere positive
Jacobian \cite{AlessandriniNesi01}; a quantitative positive lower bound on
this Jacobian would give $\min_y|Z(y,\omega)|>0$ for every $\omega$ and, by
compactness, a uniform bound; critical points of solutions to constant-coefficient anisotropic
$p$-Laplace equations in the plane are isolated \cite{AlessandriniSigalotti01}.
A corresponding statement for periodic correctors of the variable-coefficient
degenerate law is, to our knowledge, not available; \cite{Wolf22} assumes
it and refers to the two-dimensional literature for $d=2$. The averaged quantity
used in Theorem~\ref{thm:geom}(b) is more tolerant than \eqref{eq:test}: if
$Z$ vanishes only at isolated points and only to first order, then
$|\{y:|Z(y,\omega)|<t\}|\lesssim t^m$ and $\Lambda_0<\infty$ as soon as
$p<m+2$, that is $p<4$ on surfaces.
\end{remark}

\subsection{Layered materials}
For laminates, \eqref{eq:kappa} for the continuous correctors is shown in
\cite{Wolf22}. The following version gives an explicit constant, uniform
in $\delta$, and covers the discrete correctors.

\begin{proposition}[Layered materials]\label{prop:laminate}
Let $m\ge1$ and $a(y,z)=c(y_1)\big(\delta^2+c(y_1)z^\top\bar Az\big)^{\frac{p-2}{2}}\bar Az$,
that is \eqref{eq:model} with $A(y)=c(y_1)\bar A$, where $\bar A$ is a constant
symmetric positive definite matrix with condition number $\mathrm{cond}(\bar A)$
and $c$ is measurable and periodic with $0<c_-\le c\le c_+$, $K=c_+/c_-$.
Then the continuous corrected gradients satisfy \eqref{eq:kappa} with
$\kappa\ge K^{-p/2}\,\mathrm{cond}(\bar A)^{-1/2}$, for every $\delta\ge0$.
The same holds for $Z_H$ when $c$ is constant on the columns
$\{i/M<y_1<(i+1)/M\}$ of a periodic simplicial mesh that resolves these
columns, and $W_H$ is the corresponding $P_1$ space.
\end{proposition}

\begin{proof}
\emph{Continuous case.} By uniqueness the corrector depends on $y_1$ only,
$Z=\xi+t(y_1)e_1$, so $Z_k=\xi_k$ for $k\ge2$. Let
$\psi:=e_1^\top\bar AZ=\bar A_{11}Z_1+\sum_{k\ge2}\bar A_{1k}\xi_k$, an increasing
affine function of $Z_1$ with $\int_Y\psi=e_1^\top\bar A\xi$. Completing the
square,
\begin{equation}\label{eq:schur}
z^\top\bar Az=\frac{(e_1^\top\bar Az)^2}{\bar A_{11}}+z'^\top Sz',\qquad z'=(z_2,\dots,z_m),
\end{equation}
with $S$ the (positive semidefinite) Schur complement; since $Z'=\xi'$ is
constant, $X:=Z^\top\bar AZ$ is an increasing function of $|\psi|$. The cell
equation states that $\sigma:=e_1\cdot a(y,Z(y))=g(y)\psi(y)$ is constant,
where $g=c(\delta^2+cX)^{(p-2)/2}>0$. If $\sigma=0$ then $\psi\equiv0$ and
\eqref{eq:schur} gives $Z^\top\bar AZ=\xi^\top\bar A\xi$. Otherwise $\psi$ has the
sign of $\sigma$, so $|e_1^\top\bar A\xi|=\int_Y|\psi|\le\operatorname{ess\,sup}|\psi|$.
Let $y,y^*$ be points with $|\psi(y^*)|<|\psi(y)|$, and write $c=c(y)$,
$c^*=c(y^*)$, $X=X(y)$, $X^*=X(y^*)$; then $X^*<X$. Since
$g^*|\psi^*|=g|\psi|$, we need $g^*>g$, which forces $c^*>c$. Then
$\delta^2+c^*X^*\le\frac{c^*}{c}(\delta^2+cX)$ and
\[
\frac{|\psi(y)|}{|\psi(y^*)|}=\frac{g^*}{g}=\frac{c^*}{c}\Big(\frac{\delta^2+c^*X^*}{\delta^2+cX}\Big)^{\frac{p-2}{2}}\le\Big(\frac{c^*}{c}\Big)^{p/2}\le K^{p/2}.
\]
Hence $\operatorname{ess\,inf}|\psi|\ge K^{-p/2}|e_1^\top\bar A\xi|$, and by
\eqref{eq:schur}, $Z^\top\bar AZ\ge K^{-p}\xi^\top\bar A\xi$ a.e., so
$|Z|^2\ge K^{-p}\,\mathrm{cond}(\bar A)^{-1}|\xi|^2$.

\emph{Discrete case.} Let $t_H$ be the discrete minimizer among
column-wise affine functions of $y_1$; these belong to $W_H$ because the
columns are resolved. Its gradient is constant on each column, so the flux
$a(Z_H)$ is constant on each column, and the Euler equation for
perturbations by column-wise affine functions of $y_1$ makes
$e_1\cdot a(Z_H)=\sigma$ the same on all columns. For a nodal basis
function $\varphi$ of $W_H$, $\int_{\rm col}\partial_{y_2}\varphi=0$ by
periodicity in $y_2$, and $\sum_{\rm col}\sigma\int_{\rm col}\partial_{y_1}\varphi=
\sigma\int_Y\partial_{y_1}\varphi=0$. Hence $t_H$ solves \eqref{eq:cell} on
all of $W_H$, and by uniqueness it is the discrete corrector. The
continuous argument now applies column by column.
\end{proof}

The proposition covers in particular the smooth anisotropic laminate of
Section~\ref{sec:num-cell}. The exponent $p/2$ is attained only in the worst
case. For $p=2$ and $\bar A=I$ the bound $K^{-1}$ is the familiar
harmonic-mean bound, and in the regime
$|\xi|\gg\delta$, where the flux along $e_1$ behaves like $c^{p/2}|Z|^{p-2}Z_1$, the
ratio along $e_1$ behaves like $2/(K^{p/(2(p-1))}+1)$, that is $2/(K^{2/3}+1)$ for $p=4$. The
numerical experiments of Section~\ref{sec:num} reproduce both regimes.

\begin{remark}[Three dimensions]
In three dimensions the Jacobian of the linear corrector map may change sign
\cite{BrianeMiltonNesi04}, so corrected gradients can vanish at isolated
points for particular directions. By Remark~\ref{rem:open}, the averaged
index $\Lambda_0$ still remains finite at such points when $p<5$.
\end{remark}

\section{Geometric consistency on approximate surfaces}\label{sec:geometry}

\subsection{Surfaces, lifts and the two discrete problems}
Let $\Gamma\subset\mathbb R^{m+1}$ be a connected closed oriented $C^3$
hypersurface with unit normal $n$, tangential projector $P=I-n\otimes n$,
signed distance $d$ and closest-point map $p(x)=x-d(x)n(p(x))$ on a tubular
neighbourhood. Let $\Gamma_h$ be a shape-regular, quasi-uniform polyhedral
approximation with nodes on $\Gamma$ and mesh size $h$, contained in the
tubular neighbourhood, such that $p:\Gamma_h\to\Gamma$ is a bijection. For
$v:\Gamma_h\to\mathbb R$ the lift is $v^\ell(p(x))=v(x)$, and
\[
\nabla_\Gamma v^\ell(p(x))=T(x)\nabla_{\Gamma_h}v(x),\qquad
\int_\Gamma g\,d\sigma=\int_{\Gamma_h}\mu\,(g\circ p)\,d\sigma_h ,
\]
where $T=B^{-T}$, $B=P(I-dD^2d)P_h$, and $\mu$ is the area ratio.

\begin{assumption}[Geometric approximation]\label{ass:geom}
There is $C_g$ independent of $h$ such that
$\|1-\mu\|_{L^\infty(\Gamma_h)}+\|T-P\|_{L^\infty(T\Gamma_h\to T\Gamma)}\le C_gh_g^2$,
and $T$, $T^{-1}$, $\mu^{\pm1}$ are uniformly bounded.
\end{assumption}

These are the classical estimates for piecewise affine interpolation of a
$C^3$ surface \cite{Dziuk88,DziukElliott13}. For the bound on $T-P$: for
$\tau\in T\Gamma_h$, $B\tau=P\tau-dPD^2d\,\tau$ with $|d|\le Ch_g^2$, and
$L:=P|_{T\Gamma_h}:T\Gamma_h\to T\Gamma$ has singular values $1$ and
$n\cdot n_h=1-O(h_g^2)$, so $L^{-T}-L=L^{-T}(I-L^\top L)=O(h_g^2)$; we write $h_g$ for the
geometric mesh size to keep track of its role, and in our setting
$h_g\simeq h$. Table~\ref{tab:geominputs} confirms both orders for the torus
used in Section~\ref{sec:num}.

\begin{assumption}[Intrinsic discrete law]\label{ass:intrinsic}
For each $x\in\Gamma$, $q_H(x,\cdot):T_x\Gamma\to T_x\Gamma$ is the discrete
homogenized flux of a cell problem satisfying Assumption~\ref{ass:micro},
written in an orthonormal frame of $T_x\Gamma$, and the frames are such that
$q_H$ is a well-defined field on the tangent bundle. The structure constants
are uniform in $x$.
\end{assumption}

Assumption~\ref{ass:intrinsic} holds for isotropic laws in any frame, and
for anisotropic materials prescribed in a global orthonormal frame, as on
the torus of Section~\ref{sec:num}. It is restrictive for anisotropic
periodic microstructures on closed orientable surfaces. Sketch: if the
local frames carrying a fixed lattice can be chosen so that all transition
maps are lattice symmetries, the structure group of the tangent bundle
reduces to a finite subgroup $G\subset SO(2)$; the associated principal
$G$-bundle defines a finite covering $\tilde\Gamma\to\Gamma$ of degree
$k\le|G|$ on which the tangent bundle is trivial, so
$0=\chi(\tilde\Gamma)=k\,\chi(\Gamma)$ and $\Gamma$ is a torus. On surfaces
of other genus such microstructures require transition regions, which are
outside the scope of this paper.

Let $V_h\subset C^0(\Gamma_h)$ be the continuous piecewise linear space
on $\Gamma_h$ and $V_h^0=\{v\in V_h:\int_{\Gamma_h}\mu v\,d\sigma_h=0\}$, so that
$v\in V_h^0$ if and only if $v^\ell$ has mean zero on $\Gamma$. Let
$f\in L^{p'}(\Gamma)$ have mean zero and $f^e=f\circ p$. The two discrete
problems are: find $u^{\rm ex},u^{g}\in V_h^0$ such that, for all
$v\in V_h^0$,
\begin{align}
\mathcal A^{\rm ex}(u^{\rm ex};v)&:=\int_\Gamma q_H\big(x,\nabla_\Gamma(u^{\rm ex})^\ell\big)\cdot\nabla_\Gamma v^\ell\,d\sigma=\int_\Gamma f\,v^\ell\,d\sigma,\label{eq:ex}\\
\mathcal A^{g}(u^{g};v)&:=\int_{\Gamma_h} q_H\big(p(x),P\nabla_{\Gamma_h}u^{g}\big)\cdot P\nabla_{\Gamma_h}v\,d\sigma_h=\int_{\Gamma_h} f^e\,v\,d\sigma_h .\label{eq:g}
\end{align}
Problem \eqref{eq:ex} is the cell-discrete method written exactly on
$\Gamma$; it contains no geometric error. Problem \eqref{eq:g} is the
practical method: the discrete tangential gradient is projected onto
$T_{p(x)}\Gamma$, and neither $T$ nor $\mu$ is evaluated. Both are uniquely
solvable by the Browder--Minty theorem, because
$(q_H(\xi)-q_H(\eta))\cdot(\xi-\eta)\ge c_1|\xi-\eta|^p$ (Section~\ref{sec:setting}),
$q_H$ is continuous by Theorem~\ref{thm:upper}, and $P|_{T\Gamma_h}$ is
uniformly invertible, and both satisfy the $W^{1,p}$ stability
bound
\begin{equation}\label{eq:stab}
\|\nabla_\Gamma(u^{\rm ex})^\ell\|_{L^p(\Gamma)}+\|\nabla_\Gamma(u^{g})^\ell\|_{L^p(\Gamma)}\le C_{\rm st},
\end{equation}
with $C_{\rm st}$ independent of $h$ and $H$: test each problem with its
solution, use $q_H(\xi)\cdot\xi\ge c_1|\xi|^p$, the uniform bounds of
Assumption~\ref{ass:geom}, the Poincaré inequality on $\Gamma$ and Young's
inequality.

\subsection{The geometric residual}

\begin{lemma}\label{lem:residual}
Under Assumptions~\ref{ass:micro}, \ref{ass:geom} and \ref{ass:intrinsic},
for all $w,v\in V_h^0$, with $\zeta=\nabla_{\Gamma_h}w$ and
$\tau=\nabla_{\Gamma_h}v$,
\[
\mathcal A^{\rm ex}(w;v)-\mathcal A^{g}(w;v)=\int_{\Gamma_h}\mathcal R_w\cdot T\tau\,d\sigma_h,\qquad
|\mathcal R_w|\le C\,h_g^2(\delta+|\zeta|)^{p-2}|\zeta| ,
\]
for a vector field $\mathcal R_w$ that depends on $w$ only, and
\[
\Big|\int_\Gamma f\,v^\ell-\int_{\Gamma_h}f^ev\Big|\le C\,h_g^2\int_{\Gamma_h}|f^e||v| .
\]
\end{lemma}

\begin{proof}
On $\Gamma_h$, $\mathcal A^{\rm ex}(w;v)=\int_{\Gamma_h}\mu\,q_H(p,T\zeta)\cdot T\tau$,
and
\[
\begin{aligned}
&q_H(P\zeta)\cdot P\tau-\mu\,q_H(T\zeta)\cdot T\tau\\
&\quad=(1-\mu)q_H(P\zeta)\cdot P\tau\\
&\qquad+\mu\big(q_H(P\zeta)-q_H(T\zeta)\big)\cdot P\tau
+\mu\,q_H(T\zeta)\cdot(P-T)\tau .
\end{aligned}
\]
By Theorem~\ref{thm:upper}, $|q_H(\chi)|\le C_{\rm up}(\delta+|\chi|)^{p-2}|\chi|$
and
\[|q_H(P\zeta)-q_H(T\zeta)|\le C_{\rm up}W(P\zeta,T\zeta)^{p-2}|(P-T)\zeta|.\]
Assumption~\ref{ass:geom} bounds $|1-\mu|$ and $\|P-T\|$ by $Ch_g^2$ and
$|P\zeta|,|T\zeta|$ by $C|\zeta|$. Writing $P\tau=(PT^{-1})T\tau$ and
$(P-T)\tau=(PT^{-1}-I)T\tau$, where $T^{-1}:T\Gamma\to T\Gamma_h$ and
$\|PT^{-1}-I\|\le Ch_g^2$ by Assumption~\ref{ass:geom}, each of the three
terms has the form $r\cdot T\tau$ with $r$ independent of $\tau$ and
$|r|\le Ch_g^2(\delta+|\zeta|)^{p-2}|\zeta|$; $\mathcal R_w$ is minus their sum.
The load difference is $\int_{\Gamma_h}(\mu-1)f^ev$.
\end{proof}

Because $q_H(0)=0$, the bound is homogeneous: it contains
$(\delta+|\zeta|)^{p-2}|\zeta|$ and not $(1+|\zeta|)^{p-1}$. This is what
allows the weighted argument below.

\subsection{Main result}

We fix $\delta_{\max}>0$ and consider $\delta\in[0,\delta_{\max}]$; constants
called uniform below do not depend on $h$, $H$, the cell space or $\delta$
in this range.

\begin{theorem}[Geometric consistency]\label{thm:geom}
Let Assumptions~\ref{ass:micro}, \ref{ass:geom} and \ref{ass:intrinsic}
hold, let $u^{\rm ex},u^{g}$ solve \eqref{eq:ex}, \eqref{eq:g}, and write
$a=\nabla_\Gamma(u^{g})^\ell$, $b=\nabla_\Gamma(u^{\rm ex})^\ell$ and
\[
\mathcal E_{\Gamma,H}(u^{g},u^{\rm ex})^2:=\int_\Gamma d_{x,H}(a,b)^2\,d\sigma .
\]
\begin{enumerate}[label=(\alph*),leftmargin=2em]
\item (Unconditional.) If $f\in L^{p'}(\Gamma)$, then, with a uniform constant,
\[
\|a-b\|_{L^p(\Gamma)}\le C\,h_g^{2/(p-1)},\qquad
\mathcal E_{\Gamma,H}(u^{g},u^{\rm ex})\le C\,h_g^{p/(p-1)} .
\]
\item (Optimal order, a posteriori.) Let
\[
\mathcal N_h:=\int_\Gamma\Lambda_H\big(x;a(x)\big)\,(\delta+|a(x)|)^p\,d\sigma ,
\]
\[
\mathcal N_f:=\int_\Gamma\!\int_Y\big(\delta+|Z_H(x,y,a(x))|\big)^{2-p}dy\,d\sigma ,
\]
with the convention of Definition~\ref{def:lambda}; both depend only on the
practical solution. There is a uniform constant $C$ such that:
\begin{enumerate}[label=(b\arabic*),leftmargin=2.2em]
\item if the load in \eqref{eq:g} is replaced by the exactly lifted load
$\int_\Gamma fv^\ell$, then
\[\mathcal E_{\Gamma,H}(u^{g},u^{\rm ex})\le C\,h_g^{2}\,\mathcal N_h^{1/2};\]
\item for the practical load, if $m=2$, $f\in L^2(\Gamma)$ and
$\mathcal N_f<\infty$, then
\[\mathcal E_{\Gamma,H}(u^{g},u^{\rm ex})\le C\,h_g^{2}\big(\mathcal N_h+\|f\|_{L^2}^2\mathcal N_f\big)^{1/2}.\]
\end{enumerate}
In both cases, with $\mathcal N:=\mathcal N_h$ in (b1) and
$\mathcal N:=\mathcal N_h+\|f\|_{L^2}^2\mathcal N_f$ in (b2),
$\|a-b\|_{L^p(\Gamma)}\le C\,h_g^{4/p}\mathcal N^{1/p}$ and, if $\delta>0$,
$\|a-b\|_{L^2(\Gamma)}\le C\,\delta^{-(p-2)/2}h_g^2\mathcal N^{1/2}$.
\end{enumerate}
\end{theorem}

\begin{remark}[Uniformity]\label{rem:uniform}
$\mathcal N_h\le C\Lambda_H^*$ with a uniform $C$, by \eqref{eq:stab}. Hence a
bound $\Lambda_H^*\le\bar\Lambda$ that holds uniformly for a family of cell
meshes and for $\delta\in[0,\delta_{\max}]$ (for example from
Proposition~\ref{prop:laminate}) gives (b1) with order $h_g^2$ uniformly in
that family and in $\delta$. For (b2), $\mathcal N_f\le\delta^{2-p}|\Gamma|$ gives
a bound uniform in $h$ and $H$ for each fixed $\delta>0$ but not as
$\delta\to0$; for $\delta=0$ the finiteness of $\mathcal N_f$ is an additional
hypothesis, which can fail where the macroscopic gradient vanishes on a set
of positive measure.
\end{remark}

\begin{proof}[Proof of Theorem~\ref{thm:geom}]
Let $e=u^{g}-u^{\rm ex}\in V_h^0$, so $\nabla_\Gamma e^\ell=a-b$. Testing both
problems with $e$,
\begin{equation}\label{eq:X}
X:=\int_\Gamma\big(q_H(a)-q_H(b)\big)\cdot(a-b)\,d\sigma
=\mathcal A^{\rm ex}(u^{g};e)-\mathcal A^{g}(u^{g};e)
+\ell_h(e),
\end{equation}
where $\ell_h(e)=\int_{\Gamma_h}f^ee-\int_\Gamma fe^\ell$ for the practical load
and $\ell_h=0$ for the exactly lifted load. By \eqref{eq:dxH-equiv},
$X\ge c\,\mathcal E_{\Gamma,H}^2$, and by \eqref{eq:basic}, Jensen and
\eqref{eq:mono}, $X\ge c_1\|a-b\|_{L^p}^p$. By Assumption~\ref{ass:geom},
$|\nabla_{\Gamma_h}u^g|\le C|a\circ p|$.

\emph{(a)} By Lemma~\ref{lem:residual}, Hölder's inequality, \eqref{eq:stab}
and the Poincaré inequality for $e^\ell$,
$X\le Ch_g^2\big(\|\delta+|a|\|_{L^p}^{p-1}+\|f\|_{L^{p'}}\big)\|a-b\|_{L^p}\le Ch_g^2\|a-b\|_{L^p}$.
With $X\ge c_1\|a-b\|_{L^p}^p$ this gives $\|a-b\|_{L^p}\le Ch_g^{2/(p-1)}$, and
then $\mathcal E_{\Gamma,H}^2\le CX\le Ch_g^{2+2/(p-1)}$.

\emph{(b)} By Lemma~\ref{lem:residual} with $w=u^g$, $v=e$, the geometric
residual equals $\int_{\Gamma_h}\mathcal R\cdot(a-b)\circ p$, with
$|\mathcal R|\le Ch_g^2(\delta+|a\circ p|)^{p-1}$. At a fixed $x\in\Gamma$ write
$Z_a=Z_H(x,\cdot,a(x))$, $Z_b=Z_H(x,\cdot,b(x))$, $D=Z_a-Z_b$ and
$X(x)=(q_H(a)-q_H(b))\cdot(a-b)$. Since $a-b=\int_YD$, the Cauchy--Schwarz
inequality, $W(Z_a,Z_b)\ge\delta+|Z_a|$, \eqref{eq:mono} and \eqref{eq:basic} give
\begin{equation}\label{eq:L1micro}
\begin{aligned}
|a-b|&\le\int_Y|D|\\
&\le\Big(\int_Y(\delta+|Z_a|)^{2-p}\Big)^{1/2}
\Big(\int_YW(Z_a,Z_b)^{p-2}|D|^2\Big)^{1/2}\\
&\le\Big(\int_Y(\delta+|Z_a|)^{2-p}\Big)^{1/2}\big(c_1^{-1}X(x)\big)^{1/2}.
\end{aligned}
\end{equation}
Hence $|\mathcal R||a-b|\le Ch_g^2\big(\Lambda_H(x;a)(\delta+|a|)^p\big)^{1/2}X(x)^{1/2}$;
where $\delta+|a(x)|=0$ the left-hand side vanishes, consistently with the
convention. After integration, the bounded change of measure between
$\Gamma_h$ and $\Gamma$, and the Cauchy--Schwarz inequality, the geometric
residual is bounded by $Ch_g^2\mathcal N_h^{1/2}X^{1/2}$. In case (b1) this
gives $X\le Ch_g^4\mathcal N_h$. In case (b2), by Lemma~\ref{lem:residual},
the Sobolev--Poincaré inequality
$\|e^\ell\|_{L^2(\Gamma)}\le C\|\nabla_\Gamma e^\ell\|_{L^1(\Gamma)}$ for
mean-zero functions on a closed surface ($m=2$), and \eqref{eq:L1micro}
integrated over $\Gamma$,
\[
|\ell_h(e)|\le Ch_g^2\|f\|_{L^2}\|a-b\|_{L^1(\Gamma)}
\le Ch_g^2\|f\|_{L^2}\,\mathcal N_f^{1/2}X^{1/2},
\]
so $X\le Ch_g^4\mathcal N$. The bounds follow from $X\ge c\,\mathcal E_{\Gamma,H}^2$,
$X\ge c_1\|a-b\|_{L^p}^p$ and, squaring \eqref{eq:L1micro} with
$\delta+|Z_a|\ge\delta$, $\|a-b\|_{L^2}^2\le\delta^{2-p}c_1^{-1}X$. All constants
depend only on $p$, $c_1$, $C_1$, $\delta_{\max}$, Assumption~\ref{ass:geom}
and $C_{\rm st}$.
\end{proof}

\begin{remark}[What is and is not assumed]\label{rem:aposteriori}
Part (b) needs no bound on the exactly lifted solution: $\mathcal N_h$ and
$\mathcal N_f$ involve only the practical solution $u^g$ and the corrected
gradients $Z_H(x_j,\cdot,a(x_j))$ that FE--HMM computes at its quadrature
points anyway, so they are available without additional cell solves. They
are averages, not suprema: the index enters through a
$(\delta+|a|)^p$-weighted mean over $\Gamma$ of the negative moment of order
$p-2$ of the cell fields, so isolated macroscopic points or small cell
regions with small corrected gradients are tolerated. No bound on
$\|\nabla u\|_{L^\infty}$, no higher integrability of discrete gradients and
no regularity of the corrector map with respect to $\xi$ are used. For
$p=2$, $\Lambda_H\equiv1$ and part~(b1) reduces to the classical
second-order geometric estimate \cite{Dziuk88,Demlow09}.
\end{remark}

\begin{remark}[Energy argument versus operator argument]
Comparing the values of the two discrete energies, as in the variational
crime analysis of surface finite elements, perturbs the squared natural
distance by $O(h_g^2)$ and therefore gives only
$\mathcal E_{\Gamma,H}=O(h_g)$. Part~(a) improves this to
$O(h_g^{p/(p-1)})$ using operator monotonicity only. Part~(b) shows that the
remaining gap to $O(h_g^2)$ is controlled by a negative moment of the
corrected gradients, which is a sufficient condition for the lower half of
the shifted structure of the discrete homogenized law
(Proposition~\ref{prop:reduction}).
\end{remark}

\subsection{Quadrature}
In practice both forms in \eqref{eq:ex}--\eqref{eq:g} and the load are
evaluated by a quadrature rule $Q$ on the elements of $\Gamma_h$, the cell
problems are solved only at its points, and the mean constraint is imposed
through $\ell_Q(v):=\sum_jw_j\mu(x_j)v(x_j)$. We write $\int^Q$ for the
corresponding weighted sums.

\begin{lemma}[Discrete inequalities]\label{lem:quad}
Let $Q$ have positive weights, points in the element interiors, and be
exact for polynomials of degree two on each flat element (for instance the
three-point rule with barycentric points $(\frac23,\frac16,\frac16)$). Then,
uniformly in $h$, for all $v\in V_h$ with $\ell_Q(v)=0$ and $1\le r<\infty$:
(i) $\int^Q_{\Gamma_h}|\nabla_{\Gamma_h}v|^r=\int_{\Gamma_h}|\nabla_{\Gamma_h}v|^r$ and
$\int^Q_{\Gamma_h}v^2=\int_{\Gamma_h}v^2$;
(ii) $\|v\|_{L^p(\Gamma_h)}\le C\|\nabla_{\Gamma_h}v\|_{L^p(\Gamma_h)}$;
(iii) for $m=2$, $\|v\|_{L^2(\Gamma_h)}\le C\|\nabla_{\Gamma_h}v\|_{L^1(\Gamma_h)}$.
\end{lemma}

\begin{proof}
(i) holds because $\nabla_{\Gamma_h}v$ is elementwise constant, the weights on
each element sum to its area, and $v^2$ is quadratic on each element. For
(ii)--(iii), let $\bar v$ be the mean of $v$ over $\Gamma_h$. By equivalence of
norms on the reference element and shape regularity,
$|v(x_j)|\le C|T|^{-1}\|v\|_{L^1(T)}$ for points $x_j\in T$, so
$|\ell_Q(w)|\le C\|w\|_{L^1(\Gamma_h)}$ for $w\in V_h$, with $C$ uniform, and
$\ell_Q(1)\ge c>0$. From $\ell_Q(v)=0$,
$|\bar v|\,\ell_Q(1)=|\ell_Q(v-\bar v)|\le C\|v-\bar v\|_{L^1}$, and
$\|v\|\le\|v-\bar v\|+C|\bar v|$ in either norm. The mean-zero Poincaré and
Sobolev--Poincaré inequalities on $\Gamma_h$, uniform for a family of
polyhedral surfaces converging to $\Gamma$ (by transport to $\Gamma$ with
Assumption~\ref{ass:geom}), conclude.
\end{proof}

For quadrature-based loads, the data are assumed to have well-defined
point values at the quadrature nodes and a uniformly bounded discrete data
norm
\[
\|f\|_{r,Q}:=\Big(\sum_jw_j\,\mu(x_j)\,|f^e(x_j)|^r\Big)^{1/r},\qquad r\in\{2,p'\},
\]
which holds, for instance, for $f\in C^0(\Gamma)$; for $f$ only in
$L^2(\Gamma)$ the point values are not defined and a stable approximation of
the load must be used instead. With Lemma~\ref{lem:quad}, the proof of
Theorem~\ref{thm:geom} carries over to the quadrature forms with $\int$
replaced by $\int^Q$ in $X$, $\mathcal E_{\Gamma,H}$, $\mathcal N_h$ and
$\mathcal N_f$, and with $\|f\|_{L^{p'}(\Gamma)}$ in (a) and
$\|f\|_{L^2(\Gamma)}$ in (b2) replaced by $\|f\|_{p',Q}$ and $\|f\|_{2,Q}$: the pointwise
estimates \eqref{eq:L1micro} and those of Lemma~\ref{lem:residual} hold at
each quadrature point, Hölder's and the Cauchy--Schwarz inequality hold for
positive weights, the lifted gradient is a bounded perturbation of the
elementwise constant $\nabla_{\Gamma_h}v$, and (ii)--(iii) replace the
continuous Poincaré and Sobolev--Poincaré inequalities. If both problems
use the same rule, as in our experiments, the geometric error is
isolated. The difference between quadrature-based and exactly integrated
methods is a separate macroscopic quadrature error, which enters the total
estimate as an additional term and is treated as in
\cite{AbdulleVilmart14,AbdulleHuber16}; Corollary~\ref{cor:total} is stated
for exact macroscopic integration.

\subsection{Complete error estimate}

For the variational case $a=D_z\phi$ let
\[
\Phi_{\rm hom}(x,\xi)=\int_Y\phi\big(x,y,Z(x,y,\xi)\big)\,dy,\qquad
\Phi_H(x,\xi)=\int_Y\phi\big(x,y,Z_H(x,y,\xi)\big)\,dy
\]
be the continuous and discrete effective energies, $u_0\in W^{1,p}(\Gamma)$ the mean-zero minimizer of
$J(v)=\int_\Gamma\Phi_{\rm hom}(x,\nabla_\Gamma v)-\int_\Gamma fv$, and $d_x$ the
corrector-induced metric \eqref{eq:dxH} built with the continuous correctors
$Z$. Write $B_\phi(z,w)=\phi(z)-\phi(w)-D\phi(w)\cdot(z-w)$; under
Assumption~\ref{ass:micro}, $B_\phi(z,w)\simeq W(z,w)^{p-2}|z-w|^2\simeq|F_\delta(z)-F_\delta(w)|^2$
\cite{DieningEttwein08,DieningKreuzer08}.

\paragraph{Sampled coefficients}
Computations evaluate the coefficient at one point per cell element. For
the model law \eqref{eq:model} let $A_H$ be the elementwise constant
coefficient with $A_H|_T=A(y_T)$ at the barycentre $y_T$, and denote by
$\Phi_H^{A_H}$, $Z_H^{A_H}$ the discrete energy and corrected gradient
computed with $A_H$. All results of Sections~\ref{sec:structure} and
\ref{sec:geometry} up to here hold for $A_H$ in place of $A$, with the same
constants, since $A_H$ satisfies the same bounds; in particular
Theorem~\ref{thm:geom} compares $u^g$ and $u^{\rm ex}$ at fixed sampled
material. What requires an argument is the connection with the continuous
material.

\begin{lemma}[Coefficient sampling]\label{lem:sampling}
Let $A\in C^2(\overline Y)$ and let the cell meshes be shape regular with
mesh size $H$. For the model law there is $C$ independent of $H$, $\delta$ and
$\xi$ such that, for all $\xi$,
\[
|\Phi_H^{A_H}(\xi)-\Phi_H(\xi)|\le CH^2(\delta+|\xi|)^p,\]
\[
\int_YB_\phi\big(Z_H^{A_H}(\xi),Z(\xi)\big)dy\le\Phi_H(\xi)-\Phi_{\rm hom}(\xi)+CH^2(\delta+|\xi|)^p .
\]
\end{lemma}

\begin{proof}
For $w\in W_H$, $\xi+\nabla w$ is constant on each element $T$, and
$y\mapsto\phi(y,z)$ has second derivatives bounded by $C\|A\|_{C^2}(\delta+|z|)^p$.
The one-point barycentric rule is exact for affine functions, so
$|\int_T\phi(y,z_T)\,dy-|T|\phi(y_T,z_T)|\le CH^2|T|(\delta+|z_T|)^p$. Summing,
the energies $E^A(w)$ and $E^{A_H}(w)$ of $w$ for the two coefficients
differ by at most $CH^2\int_Y(\delta+|\xi+\nabla w|)^p$. Applying this at the two
minimizers, whose $p$-energies are bounded by $C(\delta+|\xi|)^p$ by
Lemma~\ref{lem:stab}, gives the first bound. For the second, by the cell
equation for $Z$, $\int_YB_\phi(z,Z(\xi))=E^A(z-\xi)-\Phi_{\rm hom}(\xi)$ for every
$z=\xi+\nabla w$, $w\in W_0$; take $z=Z_H^{A_H}(\xi)$ and use
$E^A(Z_H^{A_H}-\xi)\le E^{A_H}(Z_H^{A_H}-\xi)+CH^2(\delta+|\xi|)^p=\Phi_H^{A_H}(\xi)+CH^2(\delta+|\xi|)^p$
together with the first bound.
\end{proof}

\begin{assumption}[Approximation]\label{ass:approx}
(i) There are $v_h\in V_h^0$ with
$\|F_\delta(\nabla_\Gamma v_h^\ell)-F_\delta(\nabla_\Gamma u_0)\|_{L^2(\Gamma)}\le C_Mh$
and $\|\nabla_\Gamma v_h^\ell\|_{L^p}\le C$.
(ii) For some $s>0$ and all $x,\xi$,
$\inf_{w\in W_H}\|F_\delta(\xi+\nabla_yw)-F_\delta(Z(x,\cdot,\xi))\|_{L^2(Y)}^2\le C_\mu H^{2s}(\delta+|\xi|)^p$.
\end{assumption}

Part~(i) is the natural-variable interpolation estimate for $P_1$ elements
under $F_\delta(\nabla u_0)\in H^1$ \cite{BarrettLiu94,DieningRuzicka07}, here
on the lifted surface space. Part~(ii) is expected with $s=1$ for $P_1$
cells when $F_\delta(Z(x,\cdot,\xi))\in H^1(Y)$ with the indicated growth;
Table~\ref{tab:cell} shows observations compatible with $s=1$ for the
tested states. Since $W_H$ is conforming, Assumption~\ref{ass:approx}(ii)
and the identity $\Phi_H(\xi)-\Phi_{\rm hom}(\xi)=\min_{w\in W_H}\int_YB_\phi(\xi+\nabla w,Z(\xi))$
(from the cell equation for $Z$) give
\begin{equation}\label{eq:DeltaH}
0\le\Delta_H(\xi):=\Phi_H(\xi)-\Phi_{\rm hom}(\xi)\le CH^{2s}(\delta+|\xi|)^p .
\end{equation}

\begin{corollary}[Complete estimate]\label{cor:total}
Let $a=D_z\phi$ be the model law with $A\in C^2(\overline Y)$, let
Assumptions~\ref{ass:micro}, \ref{ass:geom}, \ref{ass:intrinsic} and
\ref{ass:approx} hold, let the cell problems use the sampled coefficient
$A_H$, and let the macroscopic integrals be exact. Assume (b1) or (b2) of
Theorem~\ref{thm:geom} with $\mathcal N\le\bar{\mathcal N}$ uniformly in $h$ and
$H$. Let $\sigma:=\min\{s,1\}$.
\begin{enumerate}[label=(\roman*),leftmargin=2em]
\item (Microfield.) The reconstructed microfield
$G_h(x,y)=Z_H^{A_H}(x,y,\nabla_\Gamma(u^g)^\ell(x))$ satisfies
\[
\Big(\int_\Gamma\!\int_Y\big|F_A(G_h)-F_A(Z(x,y,\nabla_\Gamma u_0))\big|^2\Big)^{1/2}
\le C\big(h+H^\sigma+h_g^2\big).
\]
\item (Macroscopic natural distance.) If moreover $\Theta_0^*<\infty$ for the
continuous correctors (by Proposition~\ref{prop:reduction},
$\Theta_0^*\le2^{p-2}\Lambda_0^*$), then
$\|F_\delta(\nabla_\Gamma(u^g)^\ell)-F_\delta(\nabla_\Gamma u_0)\|_{L^2(\Gamma)}\le C\big(h+H^\sigma+h_g^2\big)$.
\end{enumerate}
\end{corollary}

\begin{proof}
Write $\|\cdot\|$ for the $L^2(\Gamma\times Y)$ norm of differences of
$F_A$, and $\mathcal E_\Gamma(v,w)^2=\int_\Gamma d_x(\nabla_\Gamma v,\nabla_\Gamma w)^2$.
Constants may depend on $\bar{\mathcal N}$ and on the data, not on $h,H,h_g$.

\emph{Step 1 (Bregman lifting).} For $\xi,\eta$, the cell equation for
$Z(\eta)$ tested with $\pi_\xi-\pi_\eta$ gives
\[
\Phi_{\rm hom}(\xi)-\Phi_{\rm hom}(\eta)-q(\eta)\cdot(\xi-\eta)=\int_YB_\phi\big(Z(\xi),Z(\eta)\big)dy\simeq d_x(\xi,\eta)^2 .
\]

\emph{Step 2 (energy quasi-optimality on the surface).} Let
$J_H(v)=\int_\Gamma\Phi_H^{A_H}(x,\nabla_\Gamma v^\ell)-\int_\Gamma fv^\ell$ on $V_h^0$.
Problem \eqref{eq:ex} with sampled coefficients is the Euler equation of
$J_H$, which is strictly convex, so $u^{\rm ex}$ minimizes $J_H$. The
Euler equation of $J$ at $u_0$, tested with $w-u_0$ for mean-zero
$w\in W^{1,p}(\Gamma)$, and Step~1 give
$J(w)-J(u_0)=\int_\Gamma\int_YB_\phi(Z(\nabla w),Z(\nabla u_0))\simeq\mathcal E_\Gamma(w,u_0)^2$.
By Lemma~\ref{lem:sampling}, $\Phi_{\rm hom}\le\Phi_H\le\Phi_H^{A_H}+CH^2(\delta+|\cdot|)^p$
and $\Phi_H^{A_H}\le\Phi_H+CH^2(\delta+|\cdot|)^p$. Hence, for every $v\in V_h^0$, using
\eqref{eq:stab} for $u^{\rm ex}$ and Assumption~\ref{ass:approx}(i) for $v=v_h$,
\[
J\big((u^{\rm ex})^\ell\big)\le J_H(u^{\rm ex})+CH^2\le J_H(v)+CH^2
\le J(v^\ell)+\int_\Gamma\Delta_H(\nabla_\Gamma v^\ell)+CH^2 .
\]
Subtracting $J(u_0)$, with $v=v_h$, Corollary~\ref{cor:metric} (upper bound,
which needs no index) and \eqref{eq:DeltaH},
$\mathcal E_\Gamma(u^{\rm ex},u_0)^2\le C\big(h^2+H^{2s}+H^2\big)$.

\emph{Step 3 (microfield).} With $G^{\rm ex}=Z_H^{A_H}(\cdot,\cdot,b)$,
$\hat G^{\rm ex}=Z(\cdot,\cdot,b)$, $G_0=Z(\cdot,\cdot,\nabla u_0)$ and the
equivalence of the natural variables of $A$ and $A_H$ (both $\simeq
W^{p-2}|\cdot|^2$ with uniform constants),
\[
\|G_h-G_0\|\le C\,\mathcal E_{\Gamma,H}(u^g,u^{\rm ex})+\|G^{\rm ex}-\hat G^{\rm ex}\|+\|\hat G^{\rm ex}-G_0\| .
\]
The first term is bounded by Theorem~\ref{thm:geom}(b); the second by
Lemma~\ref{lem:sampling}, \eqref{eq:DeltaH} and \eqref{eq:stab},
$\|G^{\rm ex}-\hat G^{\rm ex}\|^2\le C\int_\Gamma\big(\Delta_H(b)+H^2(\delta+|b|)^p\big)\le C(H^{2s}+H^2)$;
the third equals $\mathcal E_\Gamma(u^{\rm ex},u_0)$ up to constants. This
proves (i).

\emph{Step 4 (macroscopic distance).} With $\hat G=Z(\cdot,\cdot,a)$,
$\mathcal E_\Gamma(u^g,u_0)\simeq\|\hat G-G_0\|\le\|\hat G-G_h\|+\|G_h-G_0\|$,
and $\|\hat G-G_h\|\le C(H^{s}+H)$ as in Step~3 with $a$ in place of $b$.
Corollary~\ref{cor:metric} for the continuous correctors gives
$\|F_\delta(a)-F_\delta(\nabla u_0)\|_{L^2}^2\le C\Theta_0^*\mathcal E_\Gamma(u^g,u_0)^2$,
which proves (ii).
\end{proof}

\section{Numerical experiments}\label{sec:num}

All computations use our own Python implementation (NumPy/SciPy). Every
number below is produced by the scripts listed in
Section~\ref{sec:repro}; the tables are generated from the stored results
without manual transcription.

\subsection{Setup}\label{sec:num-setup}
\paragraph{Surface and lift}
$\Gamma$ is the torus
\[X(\alpha,\beta)=((R+r\cos\beta)\cos\alpha,(R+r\cos\beta)\sin\alpha,r\sin\beta)\]
with $R=2$, $r=0.7$. The polyhedral surface $\Gamma_h$ interpolates $X$ on a
uniform $N\times N$ parameter grid, each cell split into two triangles. The
closest-point map is explicit: for $x\in\mathbb R^3$,
$\alpha=\operatorname{atan2}(x_2,x_1)$, $\rho=(x_1^2+x_2^2)^{1/2}$,
$\beta=\operatorname{atan2}(x_3,\rho-R)$ and $p(x)=X(\alpha,\beta)$. Its
Jacobian acts on a tangent vector $t$ of $\Gamma_h$ by
$Dp\,t=(R+r\cos\beta)(\nabla\alpha\cdot t)\,e_\alpha+r(\nabla\beta\cdot t)\,e_\beta$,
with the global orthonormal frame
$e_\alpha=(-\sin\alpha,\cos\alpha,0)$,
$e_\beta=(-\sin\beta\cos\alpha,-\sin\beta\sin\alpha,\cos\beta)$.
Thus the lifted gradient $T\nabla_{\Gamma_h}v$, the area ratio $\mu$ and the
projected gradient $P\nabla_{\Gamma_h}v$ are all evaluated exactly at every
quadrature point, and problems \eqref{eq:ex} and \eqref{eq:g} are
implemented as stated. Table~\ref{tab:geominputs} confirms the second-order
geometric inputs of Assumption~\ref{ass:geom}.

\begin{table}[t]
\caption{Geometric inputs on the polygonal torus: area ratio, lifted versus
projected gradient map, and distance to $\Gamma$ at the quadrature points.
All three are $O(h_g^2)$, as required by Assumption~\ref{ass:geom}.}
\label{tab:geominputs}
\centering\small
\begin{adjustbox}{max width=\linewidth}\begin{tabular}{rcccc}
\toprule
$N$ & $h_g$ & $\|1-\mu\|_\infty$ & $\max|T-P|/|P|$ & $\|d\|_\infty$\\
\midrule
8 & 2.0665 & $4.822\cdot10^{-1}$ & $2.619\cdot10^{-1}$ & $2.044\cdot10^{-1}$\\
16 & 1.0782 & $9.290\cdot10^{-2}$ & $6.230\cdot10^{-2}$ & $5.310\cdot10^{-2}$\\
32 & 0.5455 & $2.191\cdot10^{-2}$ & $1.529\cdot10^{-2}$ & $1.340\cdot10^{-2}$\\
64 & 0.2736 & $5.400\cdot10^{-3}$ & $3.802\cdot10^{-3}$ & $3.357\cdot10^{-3}$\\
128 & 0.1369 & $1.345\cdot10^{-3}$ & $9.492\cdot10^{-4}$ & $8.398\cdot10^{-4}$\\
\midrule
EOC (last) & & 2.01 & 2.00 & 2.00\\
\bottomrule
\end{tabular}\end{adjustbox}

\end{table}

\paragraph{Material}
The microscopic law is \eqref{eq:model} with $p=4$, prescribed in the frame
$(e_\alpha,e_\beta)$, with macroscopic factor
$\gamma(\alpha,\beta)=1+0.16\cos\beta+0.09\sin\alpha\cos\beta$; since $\gamma$
multiplies the law, $q_H(x,\xi)=\gamma(x)\,q_H^{\rm cell}(\xi)$, where
$q_H^{\rm cell}$ is the discrete flux of the cell law without the factor
$\gamma$, and only $q_H^{\rm cell}$ is computed by cell solves. The cell
tensor is either the
smooth anisotropic non-separable tensor
\begin{align*}
A_{11}&=1.40+0.35\cos(2\pi y_1)+0.12\sin(2\pi(y_1+y_2)),\\
A_{22}&=0.95+0.24\sin(2\pi y_2)-0.08\cos(2\pi(y_1-y_2)),\\
A_{12}&=0.16\sin(2\pi(y_1+y_2))+0.10\cos(2\pi y_2),
\end{align*}
which satisfies $0.62I\preceq A\preceq1.92I$, or the two-phase laminate
$A=c(y_1)I$ with $c=1$ for $y_1<\frac12$ and $c=K$ otherwise. The load is
$f=\sin\alpha+0.45\cos2\beta+0.20\sin(2\alpha+\beta)$, which has zero mean
on $\Gamma$.

\begin{figure}[t]
\centering
\includegraphics[width=\textwidth]{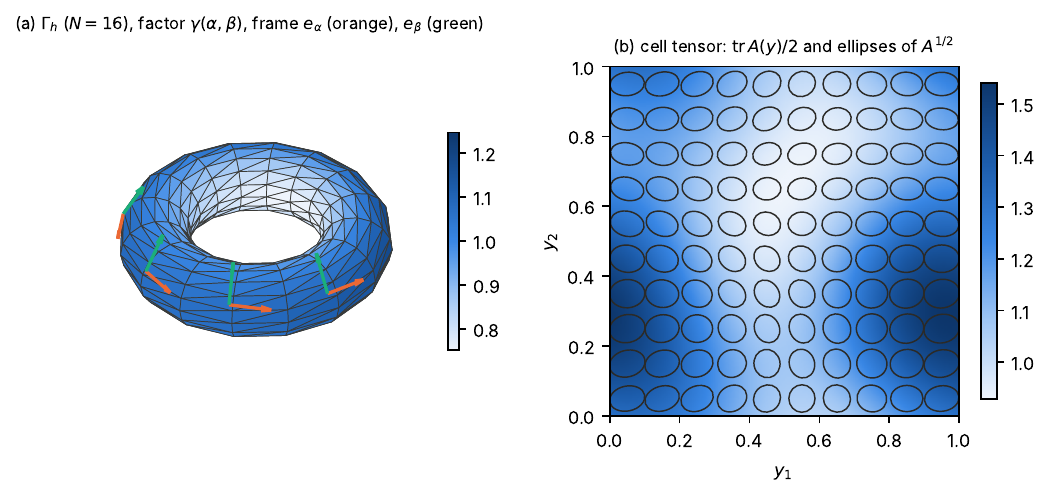}
\caption{Setup. (a) Polyhedral torus $\Gamma_h$ for $N=16$, coloured by the
macroscopic factor $\gamma(\alpha,\beta)$, with the global orthonormal frame in
which the material is prescribed. (b) Cell tensor of
Section~\ref{sec:num-setup}: half trace and ellipses of $A(y)^{1/2}$.}
\label{fig:torus-setup}
\end{figure}

Figure~\ref{fig:torus-setup} shows a coarse polyhedral surface, the
macroscopic factor, the transported frame, and the cell tensor.

\paragraph{Discretization and solvers}
Macroscopic: continuous $P_1$ on $\Gamma_h$, the $\mu$-weighted mean
constraint by a Lagrange multiplier, and the positive three-point rule of
degree two on each triangle, the same rule in both problems. Cells: periodic
$P_1$ on an $M\times M$ structured mesh ($2M^2$ triangles), the tensor
sampled at triangle barycentres, so that $q_H$ is the exact discrete
homogenized flux of the sampled material. At every quadrature point and
every macroscopic Newton iterate the cell problem is solved by a damped
Newton method from the previous iterate, and the macroscopic tangent is
assembled from the condensed cell Hessian,
$D_\xi q_H=B_{\xi\xi}-B_{\xi w}B_{ww}^{-1}B_{w\xi}$. The cell tangent agrees
with centred finite differences to relative accuracy $10^{-10}$, and the
flux agrees with the derivative of the discrete cell energy to the same
accuracy. The macroscopic Newton method uses an Armijo line search on the
discrete energy. Termination: macroscopic residual below $10^{-11}$, cell
gradient below $10^{-12}(1+|\xi|_1^3)$.

\subsection{Shifted structure of the discrete homogenized law}\label{sec:num-structure}
For the smooth material, Table~\ref{tab:structure} reports the
non-degeneracy index $\kappa_H(r)=\min|Z_H(y,\xi)|/|\xi|$, the minimum over
48 directions with $|\xi|=r$ and over all cell triangles, together with the
observed constants
\[
c_S^{\rm obs}=\min\frac{(q_H(\xi)-q_H(\eta))\cdot(\xi-\eta)}{W(\xi,\eta)^{2}|\xi-\eta|^2},\]
\[
C_{\rm up}^{\rm obs}=\max\frac{|q_H(\xi)-q_H(\eta)|}{W(\xi,\eta)^{2}|\xi-\eta|},\qquad
\varrho=\frac{d_{x,H}(\xi,\eta)^2}{W(\xi,\eta)^{2}|\xi-\eta|^2}
\]
over 400 random pairs with $|\xi|_\infty\le4$ and separations from
$10^{-3}$ to $1$. The index stays bounded away from zero over three decades
of $|\xi|$, which covers both the linear regime $|\xi|\ll\delta$ and the
quartic regime $|\xi|\gg\delta$, and both shifted constants are
essentially independent of the cell mesh, as Theorems~\ref{thm:upper} and
\ref{thm:lower} predict. The ratio $\varrho$ confirms the metric equivalence
of Corollary~\ref{cor:metric}.

\begin{table}[t]
\caption{Shifted structure of $q_H$ for the smooth anisotropic material,
$\delta=0.45$.}
\label{tab:structure}
\centering\small
\begin{adjustbox}{max width=\linewidth}\begin{tabular}{rcccccccccc}
\toprule
& \multicolumn{6}{c}{$\kappa_H$ at $|\xi|=$} & & & &\\
$M$ & 0.01 & 0.1 & 0.3 & 1 & 3 & 10 & $c_S^{\rm obs}$ & $C_{\rm up}^{\rm obs}$ & $\min\varrho$ & $\max\varrho$\\
\midrule
4 & 0.717 & 0.722 & 0.741 & 0.760 & 0.764 & 0.764 & 0.202 & 1.312 & 0.201 & 1.666\\
8 & 0.651 & 0.660 & 0.696 & 0.731 & 0.737 & 0.738 & 0.200 & 1.294 & 0.200 & 1.643\\
16 & 0.640 & 0.651 & 0.690 & 0.726 & 0.732 & 0.733 & 0.200 & 1.289 & 0.199 & 1.637\\
\bottomrule
\end{tabular}\end{adjustbox}

\end{table}

\begin{figure}[t]
\centering
\includegraphics[width=0.48\textwidth]{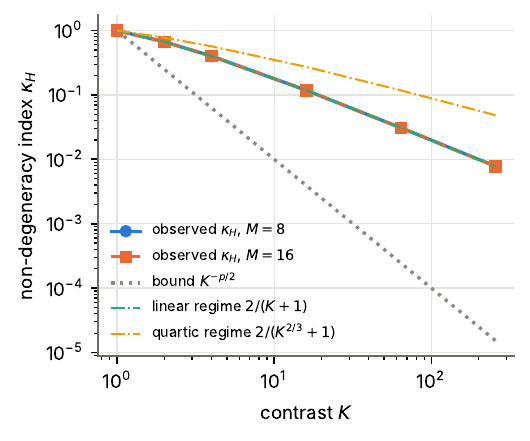}
\caption{Laminates: computed non-degeneracy index, the rigorous bound of
Proposition~\ref{prop:laminate}, and the two asymptotic regimes.}
\label{fig:laminate}
\end{figure}

For two-phase laminates, Table~\ref{tab:laminate} and
Figure~\ref{fig:laminate} compare the computed index with the bound
$K^{-p/2}=K^{-2}$ of Proposition~\ref{prop:laminate} and with the two
asymptotic regimes along $e_1$: $2/(K+1)$ for $|\xi|\ll\delta$ (harmonic
mean of $c$) and $2/(K^{2/3}+1)$ for $|\xi|\gg\delta$, where the flux along
$e_1$ behaves like $c^2|Z|^2Z_1$. The computed index reproduces both regimes
to three or four digits for $K$ from $1$ to $256$, the bound holds in every
case (it is pessimistic by a factor that approaches $2K$ for large $K$), and the values for
$M=8$ and $M=16$ coincide to rounding, as Proposition~\ref{prop:laminate}
predicts for column-aligned meshes. The lower shifted constant
$c_S^{\rm obs}$ stays of order one, and the upper constant grows like
$K^2$, as does the microscopic constant $C_1$ of this law; the normalized
ratio $C_{\rm up}^{\rm obs}/K^2$ settles at $0.34$.

\begin{table}[t]
\caption{Two-phase laminates: non-degeneracy index versus contrast.}
\label{tab:laminate}
\centering\small
\begin{adjustbox}{max width=\linewidth}\begin{tabular}{rrccccccc}
\toprule
$K$ & $M$ & $\kappa_H(0.01)$ & $2/(K+1)$ & $\kappa_H(10)$ & $2/(K^{2/3}+1)$ & bound $K^{-2}$ & $c_S^{\rm obs}$ & $C_{\rm up}^{\rm obs}/K^2$\\
\midrule
1 & 8 & 1.0000 & 1.0000 & 1.0000 & 1.0000 & $1.00\cdot10^{0}$ & 0.202 & 0.715\\
1 & 16 & 1.0000 & 1.0000 & 1.0000 & 1.0000 & $1.00\cdot10^{0}$ & 0.202 & 0.715\\
2 & 8 & 0.6669 & 0.6667 & 0.7729 & 0.7730 & $2.50\cdot10^{-1}$ & 0.363 & 0.429\\
2 & 16 & 0.6669 & 0.6667 & 0.7729 & 0.7730 & $2.50\cdot10^{-1}$ & 0.363 & 0.429\\
4 & 8 & 0.4003 & 0.4000 & 0.5681 & 0.5682 & $6.25\cdot10^{-2}$ & 0.432 & 0.364\\
4 & 16 & 0.4003 & 0.4000 & 0.5681 & 0.5682 & $6.25\cdot10^{-2}$ & 0.432 & 0.364\\
16 & 8 & 0.1178 & 0.1176 & 0.2720 & 0.2721 & $3.91\cdot10^{-3}$ & 0.530 & 0.344\\
16 & 16 & 0.1178 & 0.1176 & 0.2720 & 0.2721 & $3.91\cdot10^{-3}$ & 0.530 & 0.344\\
64 & 8 & 0.0308 & 0.0308 & 0.1176 & 0.1176 & $2.44\cdot10^{-4}$ & 1.677 & 0.342\\
64 & 16 & 0.0308 & 0.0308 & 0.1176 & 0.1176 & $2.44\cdot10^{-4}$ & 1.677 & 0.342\\
256 & 8 & 0.0078 & 0.0078 & 0.0484 & 0.0484 & $1.53\cdot10^{-5}$ & 2.850 & 0.342\\
256 & 16 & 0.0078 & 0.0078 & 0.0484 & 0.0484 & $1.53\cdot10^{-5}$ & 2.850 & 0.342\\
\bottomrule
\end{tabular}\end{adjustbox}

\end{table}

\subsection{The directional criterion for general materials}\label{sec:num-directions}
Theorem~\ref{thm:directions} reduces the uniform index bound, for a given
cell mesh, to the degenerate correctors ($\delta=0$) at unit states on a
half-circle. Table~\ref{tab:directions} reports the \emph{sampled}
directional minima $\min_k\min_T|Z_H^0|_T(\omega_k)|$, an upper bound for
$\kappa_H^0$ in \eqref{eq:test}, and the sampled maxima
$\max_k\lambda_H(0,\omega_k)$, for ten materials and $M=8,16,32$, using 90
equispaced directions $\omega_k$ for $M\le16$ and 36 for $M=32$; every degenerate cell problem is solved to a
gradient norm below $10^{-11}$. Besides the tensor of
Section~\ref{sec:num-setup}, the materials are six randomly generated
smooth anisotropic non-separable tensors
$A=R(\theta)\operatorname{diag}(\lambda_1,\lambda_2)R(\theta)^\top$, where
$\log\lambda_1$, $\log\lambda_2$ and $\theta$ are random trigonometric
polynomials of degree two scaled so that
$\lambda_i\in[10^{-1/2},10^{1/2}]$ and $\theta\in[-\pi,\pi]$, and smooth
isotropic inclusions $c(y)=1+(K-1)\exp(-(\sin^2\pi(y_1-\frac12)+\sin^2\pi(y_2-\frac12))/0.09)$
with $K=16$ (stiff), $K=1/16$ (soft) and $K=100$ (stiff). The contrast
column gives the eigenvalue ratio actually seen at the barycentres for
$M=32$.

No degeneracy is detected among the sampled directions, for any material
and resolution. These are sampled values: they detect degeneracy but do not
certify \eqref{eq:test}, for which Lemma~\ref{lem:certify} would be needed.
Between $M=16$ and $M=32$ the sampled minima change by at most $0.018$ and
the sampled maxima of $\lambda_H$ by at most $0.032$. This stability is
favourable evidence for $H$-independent limits, not a proof of uniformity
in $H$ (Remark~\ref{rem:open}). Stiff inclusions are the least favourable
case, because the field is screened inside the inclusion: $\kappa_H^0$
decreases with the contrast, to $0.094$ for $K=100$. The averaged index is
much less sensitive than the pointwise one: for $K=100$,
$(\kappa_H^0)^{-2}\approx112$ while $\sup_\omega\lambda_H(0,\omega)\approx4.4$. This is the
quantity that enters Theorem~\ref{thm:geom}(b).

\begin{table}[t]
\caption{Directional criterion of Theorem~\ref{thm:directions} for the
degenerate law ($p=4$, $\delta=0$), sampled over 90 ($M\le16$) or 36
($M=32$) directions: minimal element value of the corrected gradient and
largest negative moment over the sampled unit states. The column labelled
$\kappa_H^0$ contains sampled minima, which bound $\kappa_H^0$ from above.}
\label{tab:directions}
\centering\small
\begin{adjustbox}{max width=\linewidth}\begin{tabular}{lcccccccc}
\toprule
material & contrast & \multicolumn{3}{c}{$\kappa_H^0$} & \multicolumn{3}{c}{$\sup_\omega\lambda_H(0,\omega)$}\\
\cmidrule(lr){3-5}\cmidrule(lr){6-8}
 & & $M=8$ & $16$ & $32$ & $8$ & $16$ & $32$\\
\midrule
Sec.~5 tensor & 3.0 & 0.7371 & 0.7331 & 0.7306 & 1.048 & 1.053 & 1.054\\
random 1 & 9.9 & 0.6180 & 0.5815 & 0.5644 & 1.081 & 1.085 & 1.083\\
random 2 & 9.3 & 0.5565 & 0.4977 & 0.4898 & 1.128 & 1.131 & 1.131\\
random 3 & 8.6 & 0.6314 & 0.6101 & 0.5924 & 1.057 & 1.075 & 1.083\\
random 4 & 8.0 & 0.5325 & 0.4684 & 0.4519 & 1.120 & 1.127 & 1.126\\
random 5 & 10.0 & 0.5449 & 0.5113 & 0.5247 & 1.117 & 1.127 & 1.138\\
random 6 & 10.0 & 0.4775 & 0.4012 & 0.4122 & 1.210 & 1.220 & 1.219\\
stiff incl.\ $16$ & 15.6 & 0.3271 & 0.2746 & 0.2632 & 1.368 & 1.336 & 1.335\\
soft incl.\ $1/16$ & 11.8 & 0.8911 & 0.8064 & 0.8026 & 1.007 & 1.011 & 1.011\\
stiff incl.\ $100$ & 97.7 & 0.1205 & 0.0987 & 0.0943 & 4.629 & 4.418 & 4.386\\
\bottomrule
\end{tabular}\end{adjustbox}

\end{table}

\subsection{Cell convergence against an exact reference}\label{sec:num-cell}
For the smooth anisotropic laminate $A(y)=c(y_1)\bar A$ with
$c=1+0.6\sin(2\pi y_1)$ and $\bar A=\big(\begin{smallmatrix}1.4&0.3\\0.3&0.9\end{smallmatrix}\big)$,
the continuous corrector depends on $y_1$ only and is computed to machine
accuracy by a nested scalar Newton method on a 4096-point periodic
trapezoidal rule. Table~\ref{tab:cell} reports the natural cell error
$\|F_A(y,Z_H^{A_H})-F_A(y,Z)\|_{L^2(Y)}$, the flux error and the sampled
energy defect $\Phi_H^{A_H}-\Phi_{\rm hom}^A$ at three
tangent states. The discrete cells use the barycentric sampling $A_H$, so
the measured error includes the sampling error of
Lemma~\ref{lem:sampling}. The natural error converges with order one,
which is compatible with Assumption~\ref{ass:approx}(ii) with $s=1$ for the
tested states (the assumption itself is uniform in $x$ and $\xi$ and is not
verified by three states); flux and energy defect converge with order
two. The sampled defect is not the conforming defect $\Delta_H$ of
\eqref{eq:DeltaH}, whose sign follows from conformity; it is positive at
every level here, which is an observation, not a consequence of the
theory.

\begin{table}[t]
\caption{Cell convergence against the exact continuous laminate reference.}
\label{tab:cell}
\centering\small
\begin{adjustbox}{max width=\linewidth}\begin{tabular}{crccccccc}
\toprule
$\xi$ & $M$ & $H$ & nat.\ error & EOC & $|q_H-q|$ & EOC & $\Phi_H^{A_H}-\Phi_{\rm hom}^A$ & EOC\\
\midrule
$(0.72,-0.48)$ & 4 & 0.3536 & $2.115\cdot10^{-1}$ & -- & $4.110\cdot10^{-2}$ & -- & $8.605\cdot10^{-3}$ & --\\
 & 8 & 0.1768 & $9.495\cdot10^{-2}$ & 1.16 & $6.127\cdot10^{-3}$ & 2.75 & $1.278\cdot10^{-3}$ & 2.75\\
 & 16 & 0.0884 & $4.718\cdot10^{-2}$ & 1.01 & $1.491\cdot10^{-3}$ & 2.04 & $3.110\cdot10^{-4}$ & 2.04\\
 & 32 & 0.0442 & $2.360\cdot10^{-2}$ & 1.00 & $3.734\cdot10^{-4}$ & 2.00 & $7.786\cdot10^{-5}$ & 2.00\\
\midrule
$(0.25,0.6)$ & 4 & 0.3536 & $1.213\cdot10^{-1}$ & -- & $2.100\cdot10^{-2}$ & -- & $3.023\cdot10^{-3}$ & --\\
 & 8 & 0.1768 & $5.451\cdot10^{-2}$ & 1.15 & $3.187\cdot10^{-3}$ & 2.72 & $4.492\cdot10^{-4}$ & 2.75\\
 & 16 & 0.0884 & $2.709\cdot10^{-2}$ & 1.01 & $7.769\cdot10^{-4}$ & 2.04 & $1.093\cdot10^{-4}$ & 2.04\\
 & 32 & 0.0442 & $1.355\cdot10^{-2}$ & 1.00 & $1.945\cdot10^{-4}$ & 2.00 & $2.736\cdot10^{-5}$ & 2.00\\
\midrule
$(1.5,1)$ & 4 & 0.3536 & $1.420\cdot10^{0}$ & -- & $7.798\cdot10^{-1}$ & -- & $3.540\cdot10^{-1}$ & --\\
 & 8 & 0.1768 & $6.368\cdot10^{-1}$ & 1.16 & $1.176\cdot10^{-1}$ & 2.73 & $5.337\cdot10^{-2}$ & 2.73\\
 & 16 & 0.0884 & $3.164\cdot10^{-1}$ & 1.01 & $2.866\cdot10^{-2}$ & 2.04 & $1.300\cdot10^{-2}$ & 2.04\\
 & 32 & 0.0442 & $1.583\cdot10^{-1}$ & 1.00 & $7.174\cdot10^{-3}$ & 2.00 & $3.255\cdot10^{-3}$ & 2.00\\
\bottomrule
\end{tabular}\end{adjustbox}

\end{table}

\subsection{Geometric consistency}\label{sec:num-geom}
For each $N\in\{8,16,32,64,128\}$ we solve \eqref{eq:ex} and \eqref{eq:g} on
the same space with $M=4$, and measure $\mathcal E_{\Gamma,H}(u^g,u^{\rm ex})$
from the cell correctors at the two lifted gradients $a$ and $b$, together
with $\|a-b\|_{L^2}$, $\|a-b\|_{L^4}$ and $\|a\|_{L^6}$. Because the cells,
the quadrature and the macroscopic space are shared, the difference
isolates the geometric variational crime. The finest computation has
$98\,304$ quadrature points, each with a fully solved nonlinear cell
problem, at every macroscopic Newton iterate.

\begin{table}[t]
\caption{Geometric consistency, smooth anisotropic material, $\delta=0.45$,
$M=4$. Top: errors and experimental orders. Bottom: non-degeneracy index
and shifted monotonicity constant along the computed pair $(a,b)$, Newton
iterations, final residuals and run time.}
\label{tab:geom}
\centering\small
\begin{adjustbox}{max width=\linewidth}\begin{tabular}{rrcccccccc}
\toprule
$N$ & \#qp & $h_g$ & $\mathcal E_{\Gamma,H}$ & EOC & $\|a-b\|_{L^2}$ & EOC & $\|a-b\|_{L^4}$ & EOC & $\|a\|_{L^6}$\\
\midrule
8 & 384 & 2.0665 & $3.766\cdot10^{-1}$ & -- & $3.137\cdot10^{-1}$ & -- & $1.448\cdot10^{-1}$ & -- & 1.715\\
16 & 1536 & 1.0782 & $1.037\cdot10^{-1}$ & 1.98 & $8.209\cdot10^{-2}$ & 2.06 & $3.906\cdot10^{-2}$ & 2.01 & 1.731\\
32 & 6144 & 0.5455 & $2.707\cdot10^{-2}$ & 1.97 & $2.119\cdot10^{-2}$ & 1.99 & $1.015\cdot10^{-2}$ & 1.98 & 1.737\\
64 & 24576 & 0.2736 & $6.846\cdot10^{-3}$ & 1.99 & $5.348\cdot10^{-3}$ & 2.00 & $2.563\cdot10^{-3}$ & 1.99 & 1.739\\
128 & 98304 & 0.1369 & $1.716\cdot10^{-3}$ & 2.00 & $1.340\cdot10^{-3}$ & 2.00 & $6.423\cdot10^{-4}$ & 2.00 & 1.739\\
\bottomrule
\end{tabular}\end{adjustbox}

\vspace{0.6em}
\begin{adjustbox}{max width=\linewidth}\begin{tabular}{rcccccc}
\toprule
$N$ & $\min\kappa_H$ & $c_S^{\rm obs}$ & Newton ($g$/ex) & residual ($g$/ex) & max cell residual & time [s]\\
\midrule
8 & 0.7557 & 0.2319 & 7/8 & $2.7\cdot10^{-13}$/$7.6\cdot10^{-15}$ & $7.0\cdot10^{-14}$ & 1\\
16 & 0.7516 & 0.2316 & 8/10 & $2.3\cdot10^{-13}$/$1.5\cdot10^{-13}$ & $1.2\cdot10^{-12}$ & 7\\
32 & 0.7405 & 0.2309 & 8/8 & $2.1\cdot10^{-13}$/$2.1\cdot10^{-13}$ & $1.8\cdot10^{-12}$ & 25\\
64 & 0.7245 & 0.2309 & 8/8 & $2.4\cdot10^{-13}$/$2.4\cdot10^{-13}$ & $2.4\cdot10^{-12}$ & 154\\
128 & 0.7194 & 0.2308 & 12/8 & $7.3\cdot10^{-12}$/$2.0\cdot10^{-13}$ & $2.6\cdot10^{-12}$ & 718\\
\bottomrule
\end{tabular}\end{adjustbox}

\end{table}

\begin{figure}[t]
\centering
\includegraphics[width=\textwidth]{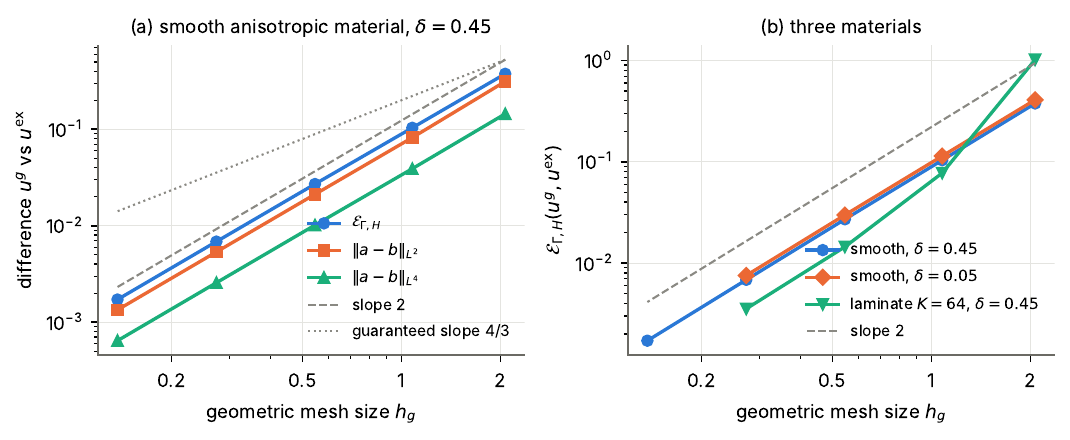}
\caption{Geometric consistency. (a) The three error measures for the smooth
material converge with order two; the guaranteed order $4/3$ of
Theorem~\ref{thm:geom}(a) is shown for comparison. (b) The corrector-induced
metric for three materials.}
\label{fig:geom}
\end{figure}

Table~\ref{tab:geom} and Figure~\ref{fig:geom} show experimental orders between $1.97$ and $2.06$
in all three measures, reaching $2.00$ on the finest level, as predicted by Theorem~\ref{thm:geom}(b). The diagnostics
explain why: the index $\kappa_H$ along the computed solution stays between
$0.72$ and $0.76$ and the monotonicity ratio near $0.23$, uniformly in $N$, so
the indices of Section~\ref{sec:structure} along the solution are bounded. The norm $\|a\|_{L^6}$ is bounded
as well, but the proof does not use it. The $L^4$ difference also converges
with order two, better than the guaranteed $h_g^{4/p}=h_g$; the theorem gives
upper bounds only.

\begin{figure}[t]
\centering
\includegraphics[width=0.92\textwidth]{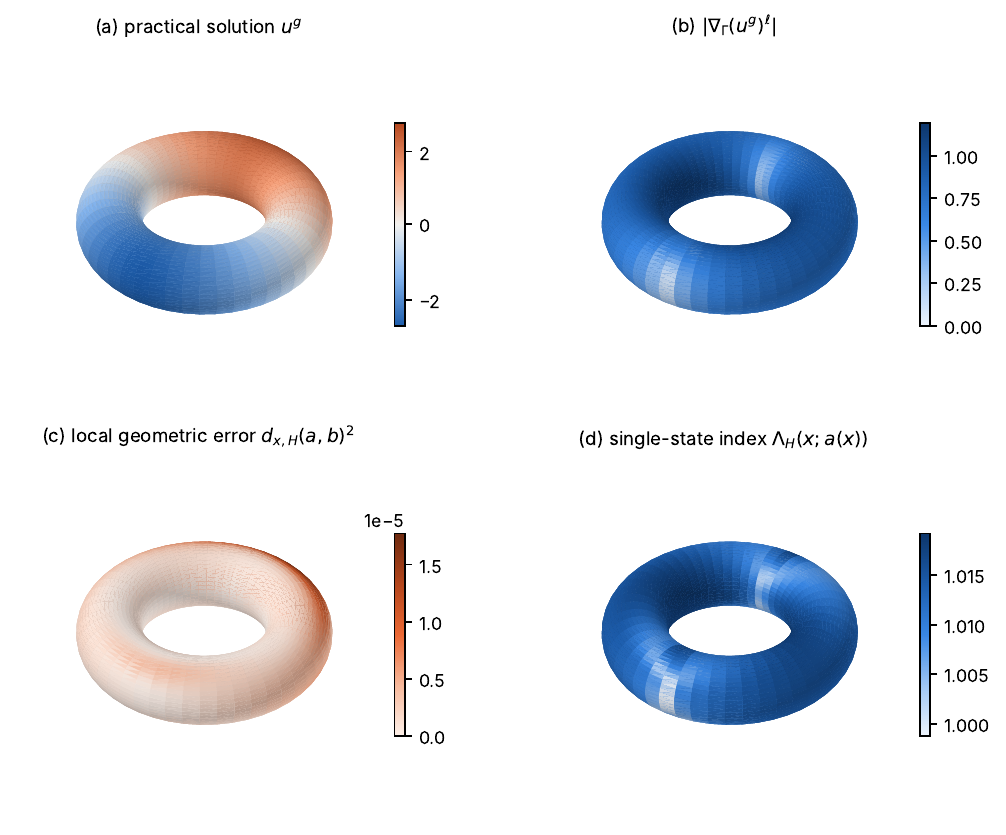}
\caption{Fields on the torus, $N=48$, $M=4$, smooth anisotropic material,
$\delta=0.45$ (element averages of quadrature values). (a) Practical solution.
(b) Lifted gradient magnitude. (c) Local density of the geometric error
$\mathcal E_{\Gamma,H}(u^g,u^{\rm ex})^2$. (d) Single-state index of
Definition~\ref{def:lambda} along the practical solution.}
\label{fig:torus-fields}
\end{figure}

\begin{figure}[t]
\centering
\includegraphics[width=\textwidth]{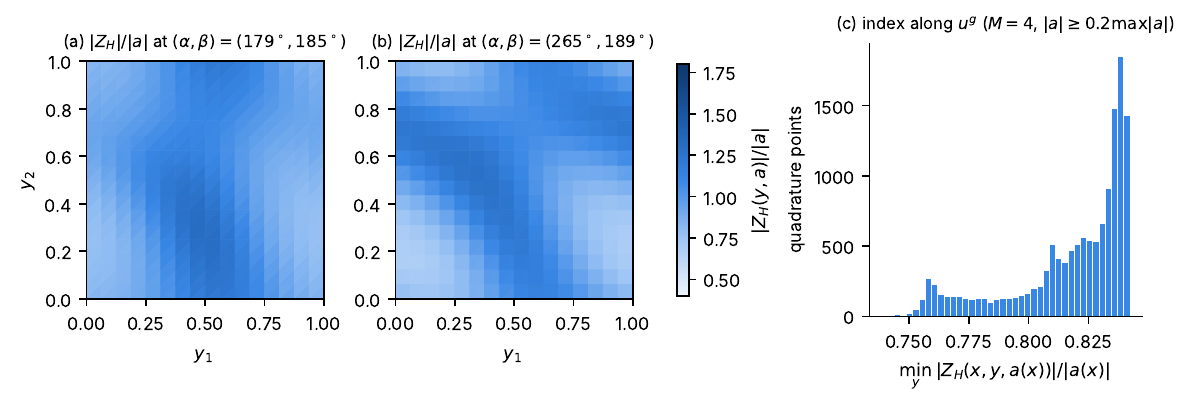}
\caption{Reconstructed microfields. (a), (b) $|Z_H(y,a(x))|/|a(x)|$ on the
cell ($M=16$) at the quadrature point of largest $|a|$ and at the point of
smallest pointwise index; $(\alpha,\beta)$ give the location on the torus.
(c) Distribution of the pointwise index $\min_y|Z_H|/|a|$ over the
quadrature points with $|a|\ge0.2\max|a|$ ($M=4$, $N=48$).}
\label{fig:torus-micro}
\end{figure}

Figure~\ref{fig:torus-fields} shows the computed fields for $N=48$
($\mathcal E_{\Gamma,H}=1.21\cdot10^{-2}$, between the values for $N=32$ and
$N=64$ in Table~\ref{tab:geom}). The local geometric error density is
spread over the whole surface, with a mean about five times larger on the
outer half ($\cos\beta>0$) than on the inner half. The single-state index
along the practical solution lies in $[0.999,1.019]$; it is closest to one
where $|a|$ is small, in the linear regime of the law, as
Proposition~\ref{prop:reduction} predicts. Figure~\ref{fig:torus-micro}
shows reconstructed microfields at two quadrature points, recomputed with
$M=16$: at the point of largest $|a|$ and at the point of smallest
pointwise index among those with $|a|\ge0.2\max|a|$, where
$\min_y|Z_H|/|a|=0.69$. Over all such quadrature points the pointwise index
lies between $0.74$ and $0.84$.

Table~\ref{tab:aposteriori} evaluates the a posteriori constants of
Theorem~\ref{thm:geom}(b) from the practical solution alone. These runs use
the practical load, so the bound to be tested is (b2), with
$\mathcal N=\mathcal N_h+\|f\|_{L^2}^2\mathcal N_f$ and $\|f\|_{L^2(\Gamma)}=5.86$; the
last column reports the corresponding ratio, the previous one the ratio
with $\mathcal N_h$ alone. For the smooth material both ratios are constant
to within $4\%$ over four refinements ($0.0071$--$0.0073$ and
$0.0023$--$0.0024$). For the laminate with $K=64$ the single-state index
along the solution stays below $3.2$, and after the pre-asymptotic first
level both ratios are bounded and decreasing.

\begin{table}[t]
\caption{A posteriori quantities of Theorem~\ref{thm:geom}(b), computed
from the practical solution only ($M=4$, $p=4$, practical load). The last
column uses the full constant $\mathcal N=\mathcal N_h+\|f\|_{L^2}^2\mathcal N_f$ of
case (b2).}
\label{tab:aposteriori}
\centering\small
\begin{adjustbox}{max width=\linewidth}\begin{tabular}{lrcccccc}
\toprule
material & $N$ & $\mathcal E_{\Gamma,H}$ & $\mathcal N_h$ & $\mathcal N_f$ & $\max_x\Lambda_H(x;a)$ & $\dfrac{\mathcal E_{\Gamma,H}}{h_g^2\mathcal N_h^{1/2}}$ & $\dfrac{\mathcal E_{\Gamma,H}}{h_g^2\mathcal N^{1/2}}$\\
\midrule
smooth, $\delta=0.45$ & 8 & $3.766\cdot10^{-1}$ & 153.04 & 37.04 & 1.019 & 0.0071 & 0.0023\\
 & 16 & $1.037\cdot10^{-1}$ & 156.27 & 37.86 & 1.019 & 0.0071 & 0.0023\\
 & 32 & $2.707\cdot10^{-2}$ & 157.55 & 38.35 & 1.019 & 0.0072 & 0.0024\\
 & 64 & $6.846\cdot10^{-3}$ & 157.95 & 38.58 & 1.019 & 0.0073 & 0.0024\\
\midrule
laminate $K=64$ & 8 & $1.017\cdot10^{0}$ & 36.54 & 133.51 & 2.359 & 0.0394 & 0.0035\\
 & 16 & $7.745\cdot10^{-2}$ & 86.43 & 126.45 & 3.002 & 0.0072 & 0.0010\\
 & 32 & $1.446\cdot10^{-2}$ & 103.52 & 127.78 & 3.119 & 0.0048 & 0.0007\\
 & 64 & $3.544\cdot10^{-3}$ & 106.95 & 128.69 & 3.130 & 0.0046 & 0.0007\\
\bottomrule
\end{tabular}\end{adjustbox}

\end{table}

Table~\ref{tab:geom-other} repeats the experiment for a nearly degenerate
law, $\delta=0.05$, and for the laminate with contrast $K=64$. Both converge
with order two. For the laminate the index along the solution lies between $0.035$
and $0.057$, more than one hundred times the worst-case bound
$K^{-2}\approx2.4\cdot10^{-4}$ of Proposition~\ref{prop:laminate}. Its coarsest levels are
pre-asymptotic (orders $3.96$ and $2.46$ in $\mathcal E_{\Gamma,H}$), after which all
three measures converge with order two. For $\delta=0.05$ and $N=64$ the exactly lifted
solve stopped at the iteration limit with residual $1.7\cdot10^{-10}$, about
seven orders of magnitude below the measured differences.

\begin{table}[t]
\caption{Geometric consistency for the smooth material with $\delta=0.05$,
with diagnostics (first two blocks) and for the laminate with $K=64$,
$\delta=0.45$, with diagnostics (last two blocks); $M=4$.}
\label{tab:geom-other}
\centering\small
\begin{adjustbox}{max width=\linewidth}\begin{tabular}{rrcccccccc}
\toprule
$N$ & \#qp & $h_g$ & $\mathcal E_{\Gamma,H}$ & EOC & $\|a-b\|_{L^2}$ & EOC & $\|a-b\|_{L^4}$ & EOC & $\|a\|_{L^6}$\\
\midrule
8 & 384 & 2.0665 & $4.088\cdot10^{-1}$ & -- & $3.401\cdot10^{-1}$ & -- & $1.592\cdot10^{-1}$ & -- & 1.817\\
16 & 1536 & 1.0782 & $1.140\cdot10^{-1}$ & 1.96 & $8.937\cdot10^{-2}$ & 2.05 & $4.287\cdot10^{-2}$ & 2.02 & 1.830\\
32 & 6144 & 0.5455 & $2.985\cdot10^{-2}$ & 1.97 & $2.307\cdot10^{-2}$ & 1.99 & $1.116\cdot10^{-2}$ & 1.98 & 1.835\\
64 & 24576 & 0.2736 & $7.556\cdot10^{-3}$ & 1.99 & $5.825\cdot10^{-3}$ & 1.99 & $2.820\cdot10^{-3}$ & 1.99 & 1.837\\
\bottomrule
\end{tabular}\end{adjustbox}

\vspace{0.5em}
\begin{adjustbox}{max width=\linewidth}\begin{tabular}{rcccccc}
\toprule
$N$ & $\min\kappa_H$ & $c_S^{\rm obs}$ & Newton ($g$/ex) & residual ($g$/ex) & max cell residual & time [s]\\
\midrule
8 & 0.7638 & 0.2843 & 6/7 & $3.1\cdot10^{-13}$/$3.8\cdot10^{-12}$ & $1.8\cdot10^{-13}$ & 4\\
16 & 0.7638 & 0.2803 & 8/8 & $7.5\cdot10^{-13}$/$6.3\cdot10^{-13}$ & $3.1\cdot10^{-12}$ & 12\\
32 & 0.7636 & 0.2641 & 9/13 & $3.9\cdot10^{-13}$/$7.0\cdot10^{-13}$ & $3.8\cdot10^{-12}$ & 60\\
64 & 0.7635 & 0.2534 & 10/39 & $3.1\cdot10^{-13}$/$1.7\cdot10^{-10}$ & $2.5\cdot10^{-12}$ & 921\\
\bottomrule
\end{tabular}\end{adjustbox}

\vspace{0.5em}
\begin{adjustbox}{max width=\linewidth}\begin{tabular}{rrcccccccc}
\toprule
$N$ & \#qp & $h_g$ & $\mathcal E_{\Gamma,H}$ & EOC & $\|a-b\|_{L^2}$ & EOC & $\|a-b\|_{L^4}$ & EOC & $\|a\|_{L^6}$\\
\midrule
8 & 384 & 2.0665 & $1.017\cdot10^{0}$ & -- & $1.554\cdot10^{-1}$ & -- & $6.187\cdot10^{-2}$ & -- & 0.795\\
16 & 1536 & 1.0782 & $7.745\cdot10^{-2}$ & 3.96 & $1.676\cdot10^{-2}$ & 3.42 & $6.738\cdot10^{-3}$ & 3.41 & 1.079\\
32 & 6144 & 0.5455 & $1.446\cdot10^{-2}$ & 2.46 & $4.292\cdot10^{-3}$ & 2.00 & $1.677\cdot10^{-3}$ & 2.04 & 1.135\\
64 & 24576 & 0.2736 & $3.544\cdot10^{-3}$ & 2.04 & $1.097\cdot10^{-3}$ & 1.98 & $4.257\cdot10^{-4}$ & 1.99 & 1.144\\
\bottomrule
\end{tabular}\end{adjustbox}

\vspace{0.5em}
\begin{adjustbox}{max width=\linewidth}\begin{tabular}{rcccccc}
\toprule
$N$ & $\min\kappa_H$ & $c_S^{\rm obs}$ & Newton ($g$/ex) & residual ($g$/ex) & max cell residual & time [s]\\
\midrule
8 & 0.0539 & 2.1034 & 8/8 & $5.0\cdot10^{-12}$/$4.7\cdot10^{-13}$ & $1.2\cdot10^{-12}$ & 2\\
16 & 0.0574 & 1.8579 & 7/7 & $9.7\cdot10^{-12}$/$8.6\cdot10^{-12}$ & $1.3\cdot10^{-12}$ & 10\\
32 & 0.0438 & 1.6809 & 8/8 & $3.3\cdot10^{-12}$/$3.3\cdot10^{-12}$ & $1.5\cdot10^{-12}$ & 50\\
64 & 0.0347 & 1.4747 & 29/8 & $2.9\cdot10^{-12}$/$2.6\cdot10^{-12}$ & $1.4\cdot10^{-12}$ & 833\\
\bottomrule
\end{tabular}\end{adjustbox}

\end{table}

\subsection{Dependence on \texorpdfstring{$\delta$}{delta} and on the load treatment}\label{sec:num-delta}
Table~\ref{tab:deltaload} repeats the geometric comparison for
$\delta\in\{0.45,0.1,0.02\}$ and $N\in\{16,32\}$ ($M=4$), with the practical load
(case (b2)) and with the exactly lifted load (case (b1)); $u^{\rm ex}$ is the
same in both cases. Between $N=16$ and $N=32$ the experimental orders of
$\mathcal E_{\Gamma,H}$ are $1.97$ for the practical load and between $2.05$ and
$2.06$ for the lifted load, for all three values of $\delta$. The ratios of
Theorem~\ref{thm:geom}(b) stay bounded as $\delta$ decreases: for (b1) the
ratio grows from $0.0068$ to $0.013$ while $\mathcal N_h$ decreases, and for
(b2) the ratio is between $0.0018$ and $0.0024$. The computed $\mathcal N_f$
grows only from $38$ to $85$; for $\delta=0.02$ this is three orders of
magnitude below the a priori bound $\delta^{2-p}|\Gamma|$ of
Remark~\ref{rem:uniform}, which is therefore very pessimistic for this
solution. To check the effect of the stopping
criterion near degeneracy we solved $\delta=0.05$, $N=32$ with macroscopic
tolerances $10^{-8}$ and $10^{-11}$: because of the quadratic convergence of
Newton's method both runs stop at the same iterate, with residuals below
$10^{-12}$, and give identical $\mathcal E_{\Gamma,H}$.

\begin{table}[t]
\caption{Geometric comparison for decreasing $\delta$ and both load
treatments ($M=4$, smooth material). The last column is the a priori bound
of $\mathcal N_f$ from Remark~\ref{rem:uniform}.}
\label{tab:deltaload}
\centering\small
\begin{adjustbox}{max width=\linewidth}\begin{tabular}{crccccccc}
\toprule
 & & \multicolumn{3}{c}{practical load (b2)} & \multicolumn{3}{c}{exactly lifted load (b1)} & \\
\cmidrule(lr){3-5}\cmidrule(lr){6-8}
$\delta$ & $N$ & $\mathcal E_{\Gamma,H}$ & $\mathcal N_f$ & $\mathcal E/(h_g^2\mathcal N^{1/2})$ & $\mathcal E_{\Gamma,H}$ & $\mathcal N_h$ & $\mathcal E/(h_g^2\mathcal N_h^{1/2})$ & $\delta^{2-p}|\Gamma|$\\
\midrule
0.45 & 16 & $1.037\cdot10^{-1}$ & 37.9 & 0.0023 & $1.034\cdot10^{-1}$ & 160.8 & 0.0070 & $2.7\cdot10^{2}$\\
0.45 & 32 & $2.707\cdot10^{-2}$ & 38.3 & 0.0024 & $2.545\cdot10^{-2}$ & 158.7 & 0.0068 & $2.7\cdot10^{2}$\\
0.1 & 16 & $1.135\cdot10^{-1}$ & 67.8 & 0.0020 & $1.041\cdot10^{-1}$ & 60.2 & 0.0115 & $5.5\cdot10^{3}$\\
0.1 & 32 & $2.973\cdot10^{-2}$ & 69.3 & 0.0020 & $2.570\cdot10^{-2}$ & 59.3 & 0.0112 & $5.5\cdot10^{3}$\\
0.02 & 16 & $1.141\cdot10^{-1}$ & 82.4 & 0.0018 & $1.042\cdot10^{-1}$ & 44.6 & 0.0134 & $1.4\cdot10^{5}$\\
0.02 & 32 & $2.988\cdot10^{-2}$ & 85.0 & 0.0018 & $2.572\cdot10^{-2}$ & 44.0 & 0.0130 & $1.4\cdot10^{5}$\\
\bottomrule
\end{tabular}\end{adjustbox}

\end{table}

\subsection{Macroscopic convergence with a manufactured solution}\label{sec:num-macro}
To check the remaining term of Corollary~\ref{cor:total}, we prescribe
$u^*=\sin\alpha+0.15\cos2\beta$ as exact solution of
$-\operatorname{div}_\Gamma(\gamma\,q_H^{\rm cell}(\nabla_\Gamma u^*))=f_H$, that is
$-\operatorname{div}_\Gamma q_H(\cdot,\nabla_\Gamma u^*)=f_H$, for the fixed cell
discretization $M=4$; the same convention is used in the code. The load $f_H$ is evaluated pointwise from the
cell flux and the cell tangent, using
$\operatorname{div}_\Gamma V=\frac{1}{R+r\cos\beta}\partial_\alpha V_1+\frac{1}{r(R+r\cos\beta)}\partial_\beta\big((R+r\cos\beta)V_2\big)$.
Table~\ref{tab:macro} reports the error of both the exactly lifted and the
practical method in the corrector-induced metric, computed from the cells at
$\nabla_\Gamma u^*$ and at the discrete gradient. Both converge with order
one, as predicted for $P_1$ elements, and their difference is far below the
macroscopic error, as Theorem~\ref{thm:geom} predicts: at practical
resolutions the geometric crime is invisible in the total error.

\begin{table}[t]
\caption{Manufactured solution for the fixed cell law, $M=4$: error of the
exactly lifted and the practical method, and the $L^2$ function error of
the practical method.}
\label{tab:macro}
\centering\small
\begin{adjustbox}{max width=\linewidth}\begin{tabular}{rcccccccc}
\toprule
$N$ & $h$ & $\mathcal E_{\Gamma,H}$ (ex) & EOC & $\mathcal E_{\Gamma,H}$ ($g$) & EOC & $\|e\|_{L^2}$ ($g$) & EOC & $|\mathcal E^{g}-\mathcal E^{\rm ex}|$\\
\midrule
8 & 2.0665 & $1.137\cdot10^{0}$ & -- & $1.208\cdot10^{0}$ & -- & $4.797\cdot10^{-1}$ & -- & $7.03\cdot10^{-2}$\\
16 & 1.0782 & $5.842\cdot10^{-1}$ & 1.02 & $5.927\cdot10^{-1}$ & 1.09 & $1.164\cdot10^{-1}$ & 2.18 & $8.57\cdot10^{-3}$\\
32 & 0.5455 & $2.942\cdot10^{-1}$ & 1.01 & $2.954\cdot10^{-1}$ & 1.02 & $3.007\cdot10^{-2}$ & 1.99 & $1.19\cdot10^{-3}$\\
64 & 0.2736 & $1.474\cdot10^{-1}$ & 1.00 & $1.475\cdot10^{-1}$ & 1.01 & $7.577\cdot10^{-3}$ & 2.00 & $1.53\cdot10^{-4}$\\
\bottomrule
\end{tabular}\end{adjustbox}

\end{table}

\subsection{Reproducibility}\label{sec:repro}
The code and numerical data are available in the public Zenodo archive at
\url{https://doi.org/10.5281/zenodo.23971110}.

The following scripts reproduce every table and figure; the raw results are
stored as JSON files and the tables are generated from them. The archived package
includes software versions, commands and random seeds.
\begin{itemize}[leftmargin=1.5em,itemsep=0pt]
\item \texttt{cells.py}: batched periodic cell solver;
\texttt{torus.py}: macroscopic solver and exact lift;
\item experiments: \texttt{e1\_geometry.py}, \texttt{e2\_delta\_load.py},
\\\texttt{n1\_structure.py}, \texttt{n3\_cell.py}, \texttt{n3\_macro.py},\\
\texttt{n4\_directions.py}, \texttt{geom\_inputs.py};
\item \texttt{make\_tables.py}, \texttt{make\_figures.py}, \texttt{fig\_torus.py}:
tables and figures.
\end{itemize}

\section{Conclusion}\label{sec:conclusion}
The geometric error of nonlinear FE--HMM on surfaces is controlled by a
property of the discrete homogenized law rather than of the macroscopic
solution. The upper half of the shifted $p$-structure is inherited, and a
negative moment of the corrected gradients is a sufficient condition for
the lower half. This moment enters the optimal second-order estimate only
as a weighted average along the practical solution, which every FE--HMM
computation evaluates without additional cell solves; for each cell mesh,
its uniform boundedness is characterized by a finite-dimensional
directional criterion, and for laminates it is explicit and uniform in the
cell mesh. Sampled tests on smooth random anisotropic materials and
inclusions of contrast up to $100$ detected no degeneracy. Open questions
are an $H$-uniform lower bound for general smooth two-dimensional
microstructures, a nonlinear and discrete counterpart of the univalence of
planar correctors, a joint refinement study in $h$ and $H$ on surfaces, and
anisotropic periodic microstructures on surfaces that do not admit a
global lattice-preserving frame.

\section*{Acknowledgements and disclosure}
Parts of the analysis, the numerical code and the text were developed with
the assistance of generative AI systems (Claude, Anthropic, and ChatGPT, OpenAI). The author
assumes responsibility for all content. The code, the raw data and the
commands needed to reproduce every table and figure are archived at \url{https://doi.org/10.5281/zenodo.23971110}.

\end{document}